\documentclass[11pt,reqno,a4paper]{amsart}
\usepackage[foot]{amsaddr}
\usepackage{a4wide}
\usepackage{amssymb}
\usepackage{graphicx}
\usepackage{subfig}
\usepackage{hyperref}
\usepackage{xcolor}
\numberwithin{equation}{section}

\advance \vsize by 1cm
\newcommand{\Nn}{\mathbb{N}}
\newcommand{\Rr}{\mathbb{R}}

\renewcommand{\epsilon}{\varepsilon}

\newcommand{\ep}{\varepsilon}

\theoremstyle{plain}
\newtheorem{theo}{Theorem}[section]
\newtheorem{lem}[theo]{Lemma}
\newtheorem{prop}[theo]{Proposition}
\newtheorem{coro}[theo]{Corollary}

\theoremstyle{definition}
\newtheorem{rem}[theo]{Remark}
\newtheorem{definition}[theo]{Definition}
\newtheorem{example}[theo]{Example}

\begin{document}

\title[A robust familly of exponential attractors]{A robust family of exponential attractors for a linear time discretization of the Cahn-Hilliard equation with a source term}

\author[D. Dor]{Dieunel Dor}
\author[M. Pierre]{Morgan Pierre}

\address{D. Dor and  M. Pierre: Laboratoire de Math\'ematiques et Applications, Universit\'e de Poitiers, CNRS, F-86073 Poitiers, France.}
\address{D. Dor: Universite d'Etat d'Haiti, Ecole Normale Superieure, Laboratoire de Mathematiques et Applications (LAMA-UEH), HT6115, 2 Rue Caseus, Pacot, Port-au-Prince, Haiti.}

\begin{abstract}We consider a linear implicit-explicit (IMEX) time discretization of the Cahn-Hilliard equation with a source term, endowed with Dirichlet boundary conditions. For every time step  small enough, we build an exponential attractor  of the discrete-in-time dynamical system associated to the discretization. We prove that, as the time step tends to 0, this attractor converges for the symmmetric Hausdorff distance to an exponential attractor of the continuous-in-time dynamical system associated with the PDE. We also prove that the fractal dimension of the exponential attractor (and consequently, of the global attractor) is bounded by a constant independent of the time step. The results also apply to the classical Cahn-Hilliard equation with Neumann  boundary conditions.
\end{abstract} 

\maketitle


\section{Introduction}
We consider the Cahn-Hilliard equation with a source term and Dirichlet boundary conditions. It reads
\begin{eqnarray}
&u_{t}+\Delta^{2}u - \Delta f(u) + g(u)=\text{ in } \Omega\times(0,+\infty),\label{a1}&\\
&u=\Delta u = 0  \text{ on }   \partial \Omega\times(0,+\infty),\label{a2}&\\
&u|_{t=0}=u_{0} \text{ in }  \Omega,\label{a3}&
\end{eqnarray}
where $\Omega$ is a bounded open subset $\mathbb{R}^{d}$ ($d=1$, $2$ or $3$) with smooth boundary $\partial\Omega$. The unknown function $u$ is the order parameter, $f$  is the nonlinear regular potential and $g$ is the source term. 

When $g=0$, the PDE~\eqref{a1} is known as the Cahn-Hilliard equation~\cite{CH} and it  has been thoroughly studied (see~\cite{Mbook} and references therein). The generalization with a source term $g$ has drawn a lot of interest in recent years, in particular for biological applications~\cite{GL17,KS08,Mbook,M21}.

 The PDE~\eqref{a1} endowed with Dirichlet boundary conditions~\eqref{a2} was analyzed in~\cite{F17,M13,M15} with various assumptions on $f$ and $g$ (see also~\cite{D22}).  In particular, global-in-time solutions were shown to exist and their asymptotic behaviour was studied. The existence of finite dimensional attractors was established.  Numerical simulations were performed, e.g., in~\cite{AKW15,F17,Letal20}.

Our purpose in this manuscript is to obtain a similar asymptotic result for a linear  time discretization of~\eqref{a1}-\eqref{a2} with fixed time step $\delta t>0$. In addition, we want a construction of exponential attractors  which is robust as  $\delta t$ goes to $0$. We use a first order implicit-explicit (IMEX) time discretization where the nonlinearities $f$ and $g$ are treated explicitly and the bilaplacian is treated implicitly. This is a very popular discretization of the classical Cahn-Hilliard equation which allows the use of the fast Fourier transform (FFT)~\cite{BPP22,SY}. It has also been successfully used in variants of the Cahn-Hilliard equation including a source term~\cite{D23,DMP}.

An exponential attractor is a compact positively invariant set which contains the global attractor, has finite fractal dimension and attracts exponentially the trajectories. In comparison with the global attractor, an exponential attractor is
expected to be more robust to perturbations: global attractors are generally upper
semicontinuous with respect to perturbations, but the lower semicontinuity can
be proved only in some specific cases (see~\cite{BV,MZ08,SH,W16} and references therein).  In particular, the upper semicontinuity of the global attractor as the mesh step and the time step tend to $0$ was proved  in~\cite{EL92} for a finite element approximation of the Cahn-Hilliard equation. 

Exponential attractors were first introduced by Eden {\it et al.}~\cite{EFNT} with a construction based on a ``squeezing property''. In~\cite{EMZ}, Efendiev, Miranville and Zelik proposed a robust construction of exponential attractors based on a ``smoothing property'' and an appropriate error estimate. Their construction has been adapted to many situations, including singular perturbations. We refer the reader to the review~\cite{MZ08}  for  details. 

In~\cite{P18}, a robust family of exponential attractors was built for a time semidiscretization of a generalized Allen-Cahn equation. An abstract result was first derived, based on the construction in~\cite{EMZ}, and it was then applied to the backward Euler scheme. The same approach was successfully applied for a time splitting scheme in~\cite{BP18},  for a discretized Ginzburg Landau equation in~\cite{B22} and for a space semidiscretization of the Allen-Cahn equation in~\cite{P18b}. In these papers, the nonlinearity was treated implicitly. Here, we also adapt the approach introduced in~\cite{P18}, but we focus on a case where the nonlinearity is treated explicitly, thus allowing a linear scheme. 

Since we use an IMEX scheme, the main condition that we impose on the potential is that $f$ is Lipschitz continous on  $\Rr$ (cf.~\eqref{a5}). This restriction can be well understood for the classical Cahn-Hilliard equation, which is a gradient flow for the $H^{-1}$ inner product, so that   there is a Lyapunov functional (the energy) naturally associated with it. In order for the IMEX scheme to have the same property, it is necessary to assume that $f$ is Lipschitz continuous  and that the time step is small enough. This is known as energy stability~\cite{BPP22,SY}. 

For $g\not=0$, the PDE~\eqref{a1}-\eqref{a2} is no longer a gradient flow and  there is no Lyapunov functional associated with it. However, the PDE is a dissipative system if $f$ satisfies a standard  dissipativity assumption (see~\eqref{a7}) and if $g$ is subordinated to $f$ (cf. Remark~\ref{remgsubo}). We prove here  that the discrete-in-time dynamical system associated to  the IMEX scheme is also dissipative if the time step is small enough (cf. Section~\ref{subsec3.2}). Typically, $f$ can be the usual cubic nonlineary which is modified into an affine function outside a compact interval as in~\eqref{FprimeK}. In turn, a typical choice for $g$ is  the symport term
\begin{equation}
\label{gsymport}
g(s)=\frac{ks}{k'+|s|}\quad (s\in\Rr),
\end{equation} 
where $k,k'>0$ ~\cite{LMG21,LMG21bis,M17}.
Our analysis also includes the case $g=0$ (the classical Cahn-Hilliard equation).

Our manuscript is organized as follows. We first give  the a priori estimates for the PDE~\eqref{a1}-\eqref{a3} in Section~\ref{sec2}. In 
Section~\ref{sec3}, we establish their discrete counterpart for the IMEX scheme. The most technical part is the dissipative $H^2$ estimate (Proposition~\ref{propH2dt}). An error estimate on finite time intervals is proved in Section~\ref{sec4}. 

The main result is given in the last section. For every time step  $ \delta t>0$ small enough, we build an exponential attractor  $ \mathcal{M}_{\delta t}$ of the discrete-in-time dynamical system associated to the IMEX scheme. We prove that $\mathcal{M}_{\delta t}$  converges to $\mathcal{M}_{0}$ for the symmmetric Hausdorff distance as $\delta t$ tends to $ 0$, where $\mathcal{M}_{0}$ is  an exponential attractor of the continuous-in-time dynamical system associated with the PDE. We also prove that the fractal dimension of $\mathcal{M}_{\delta t}$ (and consequently, of the global attractor) is bounded by a constant independent of $\delta t$. The results also apply to the classical Cahn-Hilliard equation with Neumann  boundary conditions, as pointed out in  Remark~\ref{remNeu}.

\section{The continuous problem}
\label{sec2}
\subsection{Notation and assumptions }
We make the following assumptions:
\begin{equation}\label{a4}
f \in C^{1,1}(\mathbb{R}), \ f(0)=0, 
\end{equation}
\begin{equation}\label{a5}
f'\text{ is bounded on }\Rr,
\end{equation}
\begin{equation}\label{a6}
f' \text{ is piecewise }C^1\text{ and }f''\text{ is bounded on }\Rr,
\end{equation}
\begin{equation}\label{a7}
\lim_{|s|\rightarrow \infty} \inf f'(s)>0.
\end{equation}
Assumption~\eqref{a7} is the dissipativity condition. 

The  term of symport $g$ satisfies
\begin{equation}\label{a8}
g \in C^{1}(\mathbb{R}).
\end{equation}
\begin{equation}\label{a9}
g \text{ is bounded on }\Rr,
\end{equation}
\begin{equation}\label{a10}
g' \text{ is bounded on } \mathbb{R}
\end{equation}
\begin{example}
The function $g$ defined by~\eqref{gsymport}  satisfies~\eqref{a8}-\eqref{a10}. If $g=0$,  then $g$ also satisfies~\eqref{a8}-\eqref{a10} and  equation~\eqref{a1} is the classical Cahn-Hilliard equation. 

\end{example}
\begin{rem}
We note that if $f\in C^2(\Rr)$ with $f''$  bounded on $\Rr$, then assumption~\eqref{a6} is satisfied. However, the weaker assumption~\eqref{a6} is interesting because it allows for the usual $C^1$ regularization of the cubic nonlinearity $s^3-s$, defined by
\begin{equation}
\label{FprimeK}
f(s)=f_K(s)=\begin{cases}
(3K^2-1)s-2K^3,&\quad s>K,\\
s^3-s,&\quad s\in[-K,K],\\
(3K^2-1)s+2K^3,&\quad s<-K,
\end{cases}
\end{equation}
where $K\ge 1$. 
Thus, $f_K\in C^1(\Rr)$ has a linear growth at $\pm\infty$ with 
$$\max_{s\in\Rr}|f'_K(s)|=3K^2-1.$$
This regularization is very popular for the IMEX scheme~\cite{BPP22,SY}.
\end{rem}
\begin{rem}
\label{remgsubo}
The growth of $g$ is controlled by the growth of $f$ in order to ensure that the PDE is dissipative (see~\cite[Remark 2.1]{M13}). If we suppress assumption~\eqref{a9}, then dissipativity is no longer garanteed. Indeed, let us choose $f(s)=s$ and $g(s)=-\alpha s$ ($s\in\Rr$) with $\alpha>\lambda_1^2+\lambda_1$ and $\lambda_1>0$ is the first eigenvalue of the minus Laplacian operator with Dirichlet boundary conditions. We have  $-\Delta e_1=\lambda_1 e_1$ where $e_1\in H^2(\Omega)\cap H^1_0(\Omega)$ is an eigenfunction associated to $\lambda_1$.  Then the function $u(t)=e^{\beta t}e_1$ with $\beta=\alpha-\lambda_1^2-\lambda_1>0$  solves~\eqref{a1}-\eqref{a2} but $\|u(t)\|_{L^2(\Omega)}\to+\infty$ as $t\to+\infty$.   
\end{rem}

\medskip

We set $H:=L^2(\Omega)$ and we denote by $(\cdot,\cdot)$ the scalar product both in $H$ and in $H^d$ and by $\|\cdot\|$ the induced norm. The symbol $\|\cdot\|_{X}$ will indicate the norm in the generic real Banach space $X$. Next, we set $V:=H^{1}_{0}(\Omega)$, so that $V'= H^{-1}(\Omega)$ is the topological dual of $V$. The space $V$ is endowed  with the Hilbertian norm $v\mapsto \|\nabla v\|$ which is equivalent to the usual $H^1(\Omega)$-norm, thanks to the Poincar\'e inequality.

We also denote by $A:D(A)\rightarrow{} H$ the (minus) Laplace operator $A=-\Delta$ with homogeneous Dirichlet boundary condition, with domain $D(A)=H^2(\Omega)\cap H^1_0(\Omega)$. By elliptic regularity~\cite{GT}, the norm $v\mapsto \|\Delta v\|$ is equivalent to the usual $H^2(\Omega)$-norm in $D(A)$.  Moreover, 
$$D(A^2)=\left\{v\in H^4(\Omega)\ :\ v=\Delta v=0\mbox{ on }\partial\Omega\mbox{ (in the sense of trace)}\right\},$$
and the norm $v\mapsto \|\Delta^2 v\|$ is equivalent to the usual $H^4(\Omega)$-norm in $D(A^2)$. 

It is well known that $A$ is a positive self-adjoint operator with compact resolvent, so that we can define, for $s\in \mathbb{R}$, its powers $A^{s}:D(A^{s})\rightarrow{}H $. For each $s\in\Rr$, the Hilbert space $D(A^s)$ is equipped with the norm $\|v\|_{2s}= \|A^{s}v\|$. We have $D(A^0)=H$ with $\|\cdot\|_0=\|\cdot\|$ and $D(A^{1/2})=V$ with $\|\cdot\|_1=\|\nabla \cdot\|$. Indeed, an integration by parts shows that 
\begin{equation}
\label{aux2.10}
\|\nabla v\|^2=(Av,v)=\|A^{1/2}v\|^2,\quad\forall v\in D(A).
\end{equation}

If $s_1<s_2$, then the space $D(A^{s_2})$ is continuously embedded in $D(A^{s_1})$, i.e. 
\begin{equation}
\label{sobo}
\|v\|_{2s_1}\le c_S\|v\|_{2s_2},\quad\forall v\in D(A^{s_2}),
\end{equation}
where the positive constant $c_S$ depends on $s_1$ and $s_2$. 

\medskip

By~\eqref{a4}-\eqref{a5}, the map $v\mapsto f(v) $ is Lipschitz continuous from $H$ into $H$ and from $V$ into $V$ and we have
\begin{equation}
\label{c12f}
\|f(v)-f(w)\|\le L_f\|v-w\|,\quad\forall v,w\in H,
\end{equation}
where $L_f=\sup_{s\in\Rr}|f'(s)|$. 
 Similarly, by~\eqref{a8}-\eqref{a9}, $v\mapsto g(v)$ is Lipschitz continuous from $H$ into $H$ and from $V$ into $V$ and we have
\begin{equation}
\label{c12}
\|g(v)\|\leq c_{g},\quad\forall v\in H,
\end{equation}
\begin{equation}
\label{c12bis}
\|g(v)-g(w)\|\le L_g\|v-w\|,\quad\forall v,w\in H.
\end{equation}

We define $F$ by 
\begin{equation}\label{c13}
	F(s)=\int_{0}^{s}f(t)dt
\end{equation}
We deduce from~\eqref{a7} that
\begin{equation}\label{e1}
F(s)\geq \gamma_1s^{2}-\gamma_2,\quad\forall s\in \mathbb{R}, \text{ where } \gamma_1>0,\ \gamma_2\geq 0.
\end{equation}
We deduce from~\eqref{a5}  that
\begin{equation}\label{e2}
F(s)\leq \gamma_3s^{2}+ \gamma_4,\quad\forall s\in \mathbb{R},\text{ where } \gamma_3>0,\ \gamma_4\geq 0.
\end{equation}
By~\eqref{a7}, we also have
\begin{equation}
\label{a7bis}
f(s)s\ge \gamma_5s^2-\gamma_6,\quad\forall s\in\Rr, \text{ where }\gamma_5>0,\ \gamma_6\ge 0,
\end{equation}
and 
\begin{equation}
\label{a7ter}
f'(s)\ge -\gamma_7,\quad\forall s\in\Rr,\text{ where }\gamma_7\ge 0.
\end{equation}

Assumption~\eqref{a6} implies that $f''$ has a finite number of discontinuities (the corner points of $f'$). Here, $f''$ is  the distributional derivative of $f'$ since  $C^{1,1}(\mathbb{R})=W^{2,\infty}(\Rr)$. Moreover, the following chain rule holds~\cite[Theorem 7.8]{GT}: if $v\in H^1(\Omega)$, we have $f'(v)\in H^1(\Omega)$ and  
$$\nabla f'(v)=f''(v)\nabla v.$$
Consequently, if $v\in H^2(\Omega)$, then $ f(v)\in H^2(\Omega)$ and
\begin{equation}
\label{z1}
\Delta f(v)=f''(v)\Delta v+f'(v)|\nabla v|^2.
\end{equation}
We use here that $H^1(\Omega)\subset L^4(\Omega)$ since $d\le 3$,  so that $|\nabla v|^2\in L^2(\Omega)$.

The abstract version of~\eqref{a1}-\eqref{a2} reads 
\begin{equation}
\label{a1bis}
\frac{du}{dt}+A^2u+Af(u)+g(u)=0\mbox{ in } D(A^{-1}), \mbox{ for a.e. }t>0.
\end{equation}
It is associated to the variational formulation 
$$\frac{d}{dt}(u,v)+(Au,Av)+(f(u),Av)+(g(u),v)=0\quad\mbox{ in }\mathcal{D}'(0,\infty),\ \forall v\in D(A).$$

\subsection{The continuous semigroup}
We first state the well-posedness result of our model.
\begin{theo}
\label{theoa11}
	For every $u_{0}\in V$, there exists a unique solution $u$ of~\eqref{a1}-\eqref{a3} which satisfies 
	$$u\in C^{0}([0,T] , V)\cap L^{2}(0,T; D(A))\quad\mbox{ and }\quad u_{t}\in L^{2}(0,T; V'),\quad\forall T>0.$$
	Moreover, if $u_0\in D(A)$, then $u$ satisfies
		$$u\in C^{0}([0,T], D(A))\cap L^{2}(0,T ; D(A^{2})),\quad\forall T>0.$$
\end{theo}
\begin{proof}
For $u_0\in V$, the proof of existence is based on the estimates~\eqref{a21},~\eqref{a18} and a standard Galerkin scheme. In order to prove that $u$ is continuous from $[0,T]$ into $V$, we use a  standard argument~\cite{Mbook,T}. We first show that $u$ is weakly continuous into $V$, thanks to the Strauss lemma, and  we note   that $t\mapsto \|\nabla u(t)\|^2$ is absolutely continuous since $\frac{d}{dt}\|\nabla u\|^2$ belongs to $L^1(0,T)$ by~\eqref{a19}. 
For $u_0\in D(A)$, the proof of existence is based on the estimates~\eqref{a27bis} and~\eqref{a27.4}.

 The proof of  uniqueness and of the continuous dependence with respect to the initial data in the $L^{2}$-norm follow from estimate~\eqref{a27}. 
\end{proof}

As a consequence, we have the continuous (with respect to the $L^{2}$-norm) semigroup $S_0(t)$ defined as
\begin{equation}
\label{defSt}
	S_0(t): D(A)\rightarrow D(A),\ u_0\mapsto u(t),\ t\ge 0.
\end{equation}

\subsection{Dissipative  estimates}
In this subsection, we first establish some priori estimates for the solution $u$ to the system~\eqref{a1}-\eqref{a3}. These formal estimates could be rigorously  justified by a Galerkin approximation. 
These estimates are essential in the proof of  Theorem~\ref{theoa11}.

\begin{prop}[Dissipative estimate in $L^2(\Omega)$]
 \label{propa21}
	We have 
	\begin{equation}
	\label{a21}
	\|u(t)\|^2+e^{-\epsilon_0 t}\int_0^t\|\Delta u(s)\|^2ds\le C_0\|u(0)\|^2e^{-\epsilon_0 t}+M_0,\quad t\ge 0,
	\end{equation}
	where the positive constants $\epsilon_0$, $C_0$  and $M_0$ are independent of $u(0)$. 
\end{prop}
\begin{proof}
On multiplying~\eqref{a1} by $(-\Delta)^{-1}u$ and integrating over $\Omega$, we have
\begin{eqnarray}\label{a13}
\frac{d}{dt}\|u\|_{-1}^{2} + 2\|\nabla u||^{2} + 2(f(u),u) + 2(g(u),(-\Delta)^{-1}u)=0
\end{eqnarray}
Using~\eqref{a7bis}, we find 
\begin{align*}
(f(u),u)&\geq c_{1}\int_{\Omega}u^{2}dx -c_{2}\\&\geq c_{1}\|u\|^{2} -c_{2},\quad c_{1}=\gamma_5>0,\ c_{2}\geq 0.
\end{align*}
Using \eqref{c12}, the continuous injection $H\subset D(A^{-1})$ and  Young's inequality, we find
\begin{align*}
|(g(u),(-\Delta)^{-1}u)|&\leq \|g(u)\|~\|(-\Delta)^{-1}u\|\\
&\leq c_{g}c_S\|u\|\\
&\leq \frac{c_{1}}{2}\|u\|^{2} + \frac{c_{g}^{2}c_S^2}{2c_{1}}.
\end{align*}
Combining the above estimates in~\eqref{a13}, we obtain
\begin{eqnarray}\label{a14}
\frac{d}{dt}\|u(t)\|_{-1}^{2} + c_{3}(\|\nabla u(t)\|^{2} + \|u(t)\|^{2})\leq c_{4},
\end{eqnarray}
where $ c_{3}=\min(2, c_{1})>0$ and  $c_{4}=c_{4}(c_S,c_{g}, c_{1}, c_{2})\geq 0$.

Next, we multiply~\eqref{a1} by $u$ and we integrate over $\Omega$.  We get
 \begin{eqnarray}\label{a15}
 \frac{d}{dt}\|u\|^{2} +2\|\Delta u\|^{2} +2(\nabla  f(u),\nabla u)+2(g(u),u)=0.
 \end{eqnarray}
 Thanks to~\eqref{a7ter}, we have 
 \begin{align*}
(f'(u)\nabla u ,\nabla u)&\geq -\gamma_7\|\nabla u\|^{2},\quad\gamma_7>0.
 \end{align*}
 Using \eqref{c12}, Young's inequality and the Poincaré inequality, we deduce that
 \begin{align*}
 |(g(u),u)|&\leq \|g(u)\|~\|u\|\\
 &\leq c_{g}\|u\|\\
 &\leq \gamma_7\|u\|^{2} + \frac{c_{g}^{2}}{4\gamma_7}.
 \end{align*}
 Combining the above estimates in~\eqref{a15}, we find
 \begin{eqnarray}\label{a16}
 \frac{d}{dt}\|u(t)\|^{2}+ 2\|\Delta u(t)\|^{2}\leq 2 \gamma_7(\|\nabla u(t)\|^{2}+ \|u(t)\|^{2}) + c_{5},\quad c_{5}=c_{5}(c_{g} , \gamma_7)\geq 0.
 \end{eqnarray}
 Summing~\eqref{a14} and $\alpha$ times~\eqref{a16} where $\alpha >0$ is small enough, we conclude
 \begin{eqnarray}\label{a17}
 \frac{d}{dt}E_{1}(t) + c_{6}(E_{1}(t)+ \|\Delta u(t)\|^{2})\leq c_{7},\quad c_{6}>0,\  c_{7}\geq 0,
 \end{eqnarray}
 where
 \begin{equation*}
 E_{1}(t) =\alpha\|u(t)\|^{2}+ \|u(t)\|^{2}_{-1}.
 \end{equation*}
 In particular,
 \begin{equation*}
 E_{1}(t)\geq \alpha \|u(t)\|^{2}.
 \end{equation*}
Applying Gronwall's lemma to~\eqref{a17} (see, e.g.,~\cite{T}) yields 
\begin{equation}\label{c16}
E_{1}(t)+ c_{6}e^{-c_{6} t}\int_{0}^{t}\|\Delta u(s)\|^{2}ds\leq E_{1}(0)e^{-c_{6}t}+ \frac{c_{7}}{c_{6}},\quad t\geq 0.
\end{equation}
This proves Proposition~\ref{propa21}. 
\end{proof}

\begin{prop}[Dissipative estimate in $H^{1}(\Omega)$]
\label{propa18} 
We have 
\begin{equation}
\label{a18}
\|\nabla u(t)\|^2+e^{-\epsilon_1 t}\int_0^t\|u_t(s)\|^2_{-1}ds\le C_{1}\|\nabla u(0)\|^2e^{-\epsilon_1 t}+M_1,\quad t\ge 0,
\end{equation}
where the positive constants $\epsilon_1$, $C_{1}$ and $M_1$ are independent of $u(0)$. 
\end{prop}
\begin{proof}
Testing~\eqref{a1} by $(-\Delta)^{-1}u_{t}$ and integrating over $\Omega$, we have
\begin{equation}\label{a19}
\frac{d}{dt}\|\nabla u\|^{2} + 2\|u_{t}\|^{2}_{-1}+ 2(f(u),u_{t})+ 2(g(u),(-\Delta)^{-1}u_{t})=0.
\end{equation}
Using~\eqref{c12}, the  injection $ H^{-1}\subset D(A^{-1})$ and Young's inequality, we get
\begin{align*}
|(g(u),(-\Delta)^{-1}u_{t})|&\leq \|g(u)\|~\|(-\Delta)^{-1}u_{t}\|\\
&\leq c_{g}\|u_{t}\|_{-1}\\
&\leq \frac{1}{2} \|u_{t}\|_{-1}^{2} + \frac{c_{g}^{2}}{2}.
\end{align*}
By~\eqref{c13}, we get 
\begin{align*}
(f(u),u_{t})&=\frac{d}{dt}\int_{\Omega}F(u)dx.
\end{align*}
This yields
\begin{equation}
\label{a20}
\frac{d}{dt}\left(\|\nabla u(t)\|^{2}+ 2(F(u),1)\right) + \|u_{t}(t)\|_{-1}^{2}\leq c_{g}^2,\quad t\ge 0.
\end{equation}
Adding~\eqref{a17} and $\beta$ times~\eqref{a20}, where $\beta >0$ is small enough, we deduce that
\begin{eqnarray}\label{c14}
\frac{d}{dt}E_{2}(t) + c_{8}(E_{1}(t)+\|\Delta u(t)\|^{2}+ \|u_{t}(t)\|^{2}_{-1})\leq c_{9},\quad t\ge 0,
\end{eqnarray}
where $c_{8}=\min (c_{6},\beta)>0$, $c_{9}=c_{9}(c_{7},c_{g})\geq 0$ and 
\begin{equation*}
E_{2}(t)=E_{1}(t)+\beta\|\nabla u(t)\|^{2}+ 2\beta(F(u),1).
\end{equation*}
By~\eqref{e2} and the Poincar\'e inequality, we have
\begin{equation}
\label{auxE2}
E_{2}(t)\leq c\|\nabla u(t)\|^{2}+c',\quad c,~c'>0,
\end{equation}
 and since $D(A)\subset H^1_0(\Omega)$, we deduce from~\eqref{c14} that
\begin{eqnarray}\label{c15}
\frac{d}{dt}E_{2}(t)+c_{10}(E_{2}(t)+\|u_{t}(t)\|^{2}_{-1})\leq c_{11},\quad c_{10}>0,\ c_{11}\ge 0,\quad t\ge 0.
\end{eqnarray}
Gronwall's lemma yields
\begin{equation}
\label{Gron1}
E_2(t)+c_{10}e^{-c_{10}t}\int_0^t\|u_t(s)\|^2_{-1}ds\le E_2(0)e^{-c_{10}t}+\frac{c_{11}}{c_{10}},\quad t\ge 0.
\end{equation}
By~\eqref{e1}, we have
\begin{equation}
E_{2}(t)\geq c_{12}(\|u(t)\|^{2}+\|\nabla u(t)\|^{2}) -c_{13},\quad c_{12}>0,\  c_{13}\geq 0.
\end{equation}
This estimate,~\eqref{Gron1} and~\eqref{auxE2} give  the result of Proposition~\ref{propa18}.
\end{proof}

\begin{prop}[Dissipative estimate in $H^{2}(\Omega)$]
\label{propa27bis}
We have
\begin{equation}
\label{a27bis}
\|\Delta u(t)\|^2\le Q_2(\|\Delta u(0)\|)e^{-\epsilon_2 t}+M_2,\quad t\ge 0,
\end{equation}
and
\begin{equation}
\label{a27.4}
\int_0^t\|\Delta^2 u(s)\|^2ds\le Q_2(\|\Delta u(0)\|)+C_{2}t,\quad t\ge 0,
\end{equation}
where the monotonic function $Q_2$ and  the positive constants $\epsilon_2$, $M_2$ and $C_2$ are independent of $u(0)$. 
\end{prop}
\begin{proof}
We multiply~\eqref{a1} by $\Delta^{2}u$ and integrate over $\Omega$, we have
\begin{eqnarray}\label{a28}
\frac{d}{dt}\|\Delta u\|^{2} +2\|\Delta^{2}u\|^{2} + 2(g(u),\Delta^{2}u)=2(\Delta f(u),\Delta^{2}u).
\end{eqnarray}
Using \eqref{c12} and  Young's inequality,  we find
\begin{align*}
|(g(u),\Delta^{2}u)|&\leq \|g(u)\|~\|\Delta^{2}u\|\\
&\leq c_{g}\|\Delta^{2}u\|\\
&\leq \frac{1}{4}\|\Delta^{2}u\|^{2}+ c^{2}_{g}.
\end{align*}
Moreover, we have
\begin{align*}
|(\Delta f(u),\Delta^{2}u)|&\leq \|\Delta f(u)\|~\|\Delta^{2}u\|\\&\leq \|\Delta f(u)\|^{2}+ \frac{1}{4}\|\Delta^{2}u\|^{2}.
\end{align*}	
Thus,
\begin{eqnarray}
\frac{d}{dt}\|\Delta u\|^{2} +\|\Delta^{2}u\|^{2} \leq 2\|\Delta f(u)\|^{2} + 2c^{2}_{g}.\label{a29}
\end{eqnarray}
Using the chain rule~\eqref{z1}, assumptions~\eqref{a5}-\eqref{a6} and  interpolation inequalities,  we obtain (see \cite{Mbook,T})
\begin{align*}
\|\Delta f(u)\|&\leq \|f'(u)\|_{L^\infty(\Omega)}~\|\Delta u\| + \|f''(u)\|_{L^\infty(\Omega)}~\|\nabla u\|^{2}_{L^{4}(\Omega)}\\
&\leq c_{f'}\|\Delta u\| + c_{f''}\|\nabla u\|^{2}_{L^{4}(\Omega)}\\
&\leq c'_{f'}\|u\|^{\frac{2}{3}}_{H^{1}(\Omega)}\|u\|^{\frac{1}{3}}_{H^{4}(\Omega)} + c'_{f''}\|u\|^{2}_{H^{\frac{7}{4}}},\quad H^{\frac{3}{4}} \subset L^{4}(\Omega)\\
&\leq c'_{f'}\|u\|^{\frac{2}{3}}_{H^{1}(\Omega)}\|u\|^{\frac{1}{3}}_{H^{4}(\Omega)} + c''_{f''}\|u\|^{\frac{3}{4}}_{H^{1}(\Omega)}\|u\|^{\frac{1}{4}}_{H^{4}(\Omega)}\\
&\leq (\text{Thanks~to~estimate~\eqref{a18}})\\
&\leq c''_{f'}(R_{1})\|u\|^{\frac{1}{3}}_{H^{4}(\Omega)} + c'''_{f''}(R_{1})\|u\|^{\frac{1}{4}}_{H^{4}(\Omega)},\quad t\ge 0,
\end{align*}
where $R_1=\|\nabla u(0)\|$. Since the norm $v\mapsto \|\Delta^2 v\|$ is equivalent to the usual $H^4(\Omega)$-norm in $D(A^2)$, we have
\begin{equation}
\label{a30}
\|\Delta f(u)\|^{2}\le c_{14}\|\Delta^2u\|^{\frac{2}{3}}+c_{15}\|\Delta^2u\|^{\frac{1}{2}},
\end{equation}
where $c_{14}=c_{14}(f',R_1)$ and $c_{15}=c_{15}(f'',R_1)$. 
Hence, by  Young's inequality,
\begin{align}
\|\Delta f(u)\|^{2} 
&\leq \frac{1}{8}\|\Delta^2u\|^{2}+\frac{2}{3}\left(\frac{8}{3}\right)^{\frac{1}{2}}c^{\frac{3}{2}}_{14}+ \frac{1}{8}\|\Delta^2u\|^{2} + \frac{3}{4}\times 2^{\frac{1}{3}} c^{\frac{4}{3}}_{15}\nonumber\\
&\leq \frac{1}{4}\|\Delta^2u\|^{2} + c_{16},\quad c_{16}=c_{16}(c_{14},c_{15})\geq 0.\label{a30bis}
\end{align}
This estimate, combined with~\eqref{a29}, yields
\begin{equation}\label{a31}
\frac{d}{dt}\|\Delta u(t)\|^{2} +\frac{1}{2}\|\Delta^{2}u(t)\|^{2}ds \leq  c_{17},\quad  c_{17}=2c_{g}^2+c_{16}(R_1)\geq 0,\quad t\geq 0.
\end{equation}

The dissipative estimate~\eqref{a18}  shows that there exists a time $t_1=t_1(R_1)$ such that $\|\nabla u(t)\|^2\le C_{1}+M_1$, for all $t\ge t_1$.  Integrating~\eqref{a31} on the interval $[0,t]$ for $t\le t_1$, we obtain 
\begin{equation}
\label{a31zero}
\|\Delta u(t)\|^2+\frac{1}{2}\int_0^t\|\Delta^2u(s)\|^2ds\le \|\Delta u(0)\|^2+c_{17}t_1(R_1),\quad t\in [0,t_1].
\end{equation}
For $t\ge t_1(R_1)$, the constants $c_{14}$ and $c_{15}$ in~\eqref{a30} do not depend on $R_1$ and consequently, we may recast~\eqref{a31} as 
\begin{equation}\label{a31bis}
\frac{d}{dt}\|\Delta u(t)\|^{2}+ c_{18}\|\Delta u(t)\|^{2}+\frac{1}{4}\|\Delta^{2}u(t)\|^{2} \leq  c_{17},\quad t\ge t_1(R_1),
\end{equation}
where $c_{18}>0$ and $c_{17}$ do not depend on $R_1$. 
Applying Gronwall's lemma, we obtain 
\begin{equation}
\label{a31ter}
\|\Delta u(t)\|^{2}\le \|\Delta u(t_1)\|^2 e^{-c_{18}(t-t_1)}+\frac{c_{17}}{c_{18}},\quad t\ge t_1(R_1).
\end{equation}

By integration on $[t_1,t]$, we also deduce from~\eqref{a31bis} that 
\begin{equation}
\label{a31quattro}
\frac{1}{4}\int_{t_1}^t\|\Delta^2u(s)\|^2ds\le \|\Delta u(t_1)\|^2+c_{17}(t-t_1),\quad t\ge t_1(R_1).
\end{equation}
The dissipative estimate~\eqref{a27bis} follows from~\eqref{a31zero} and~\eqref{a31ter}. Estimate~\eqref{a27.4} follows from~\eqref{a31zero}   and~\eqref{a31quattro}. 
\end{proof}

We deduce from Proposition~\ref{propa27bis} the existence of bounded absorbing set in $D(A)$ and consequently, of a global attractor associated with our semigroup $S_0(t)$~\cite{T}. 
\begin{theo}The semigroup $S_0(t)$ has a global attractor $\mathcal{A}\subset D(A)$ which is invariant ($S_0(t)\mathcal{A}=\mathcal{A}$),  bounded in $H^2(\Omega)$, compact in $L^2(\Omega)$, and which attracts the bounded sets of $D(A)$ for the $L^2(\Omega)$-norm.
\end{theo}

The following estimate shows that the semigroup is H\"older continuous in time.
\begin{lem}\label{a38}
 Let $T>0$.  If $u_{0}\in H^{2}(\Omega)\cap H^{1}_{0}(\Omega)$, then 
 \begin{equation}
 \|u(t_{1})-u(t_{2})\|^{2}\leq Q(T, \|u_{0}\|_2)|t_{1}-t_{2}|,\quad\forall t_1,t_2\in[0,T].
 \end{equation}

\end{lem}
\begin{proof}
We multiply~\eqref{a1} by $u_{t}$ and integrate over $\Omega$. We find
\begin{equation}
\label{a39}
\frac{d}{dt}\|\Delta u(t)\|^{2} + 2\|u_{t}(t)\|^{2}- 2(\Delta f(u), u_{t})+ 2(g(u),u_{t})=0.
\end{equation}
Using the Cauchy-Schwarz inequality, ~\eqref{c12} and Young's inequality, we obtain
\begin{align*}
|(g(u),u_{t})|&\leq \|g(u)\|~\|u_{t}\|\\
&\leq c_{g}\|u_{t}\|\\
&\leq \frac{1}{4}\|u_{t}\|^{2}+ c^{2}_{g}.
\end{align*}
Similarly, we have 
\begin{align*}
|(\Delta f(u),u_{t})|&\leq \|\Delta f(u)\|~\|u_{t}\|\\
&\leq \frac{1}{4}\|u_{t}\|^{2} + \|\Delta f(u)\|^{2}\\
&\leq (\text{Estimate~\eqref{a30bis}})\\
&\leq \frac{1}{4}\|u_{t}\|^{2} +\frac{1}{2}\|\Delta^2u\|^2+c_{16}(R_1),
\end{align*}
where $R_1=\|\nabla u(0)\|$. 
Combining the above estimates in~\eqref{a39}, we find
$$\frac{d}{dt}\|\Delta u(t)\|^{2} + \|u_{t}(t)\|^{2}\leq \|\Delta^2u(t)\|^2 +2c_{16}(R_1),\quad t\geq 0.$$
Integrating on $[0,T]$, we obtain
$$\|\Delta u(T)\|^2+\int_0^T\|u_t(t)\|^2dt\le \int_0^T\|\Delta^2u(t)\|^2dt+2c_{16}(R_1)T.$$
Thus, by Proposition~\ref{propa27bis},
\begin{equation}
\label{a40bis}
\int_0^T\|u_t(t)\|^2dt\le Q(\|u(0)\|_2,T),
\end{equation}
where $Q$ is a continuous and monotonic function of its arguments.

Let $t_1,t_2\in[0,T]$. Using the  Cauchy-Schwarz inequality, we have
\begin{align}
\|u(t_{1})-u(t_{2})\|&=\left\|\int_{t_{1}}^{t_{2}}u_t(s)ds\right\|\nonumber\\
&\le \left|\int_{t_{1}}^{t_{2}}\left\|u_t(s)\right\|ds\right|\nonumber\\
&\leq |t_{1}-t_{2}|^{\frac{1}{2}}\left|\int_{t_{1}}^{t_{2}}\|u_t(s)\|^{2}ds\right|^{\frac{1}{2}}.\label{41}
\end{align}
 Lemma~\ref{a38} follows from~\eqref{a40bis} and~\eqref{41}.
\end{proof}

\subsection{Estimates for the difference of two solutions}
Let now $u_{1}$ and $u_{2}$ be two solutions of system~\eqref{a1}-\eqref{a3} with initial data $u_{0,1}$ and $u_{0,2}$, respectively. 
 We set $u=u_{1}-u_{2}$ and $u_{0}=u_{0,1}-u_{0,2}$ and we have, for all $T>0$, 
 \begin{equation}\label{a22}
 u_{t}+A^{2}u  +A( f(u_{1})-f(u_{2})) + (g(u_{1})-g(u_{2}))=0\  \text{ in }\  L^2(0,T;D(A^{-1})),
 \end{equation}
 \begin{equation}\label{a23}
 u|_{t=0}=u_{0}(=u_{0,1}-u_{0,2}) \  \text{ in }\  V.
 \end{equation}

We first prove:
\begin{lem}[Uniqueness]
For all $t\ge 0$, we have
\begin{equation}\label{a27}
 \|u_{1}(t)-u_{2}(t)\|^{2} +\int_{0}^{t}\|\Delta u(s)\|^{2}ds\leq e^{c_{f,g} t}\|u_{0,1}-u_{0,2}\|^{2},
 \end{equation}
 where the positive constant $c_{f,g}$ depends only on $L_f$ and $L_g$. 
\end{lem}
\begin{proof}
On  multiplying~\eqref{a22} by $u$ in $H$, we obtain
 \begin{eqnarray}\label{a25}
 \frac{d}{dt}\|u\|^{2} +2\|Au\|^{2}+ 2(f(u_{1})-f(u_{2}),A u) + 2(g(u_{1})-g(u_{2}),u)=0.
 \end{eqnarray}
 From~\eqref{c12bis}, we deduce that
 \begin{align*}
 |(g(u_{1})-g(u_{2}),u)|&\leq \|g(u_{1})-g(u_{2})\|~\|u\|\\
 &\leq L_{g}\|u\|^{2}.
 \end{align*}
 Using~\eqref{c12f} and Young's inequality, we find
 \begin{align*}
 |(f(u_{1})-f(u_{2}),A u)|&\leq \|f(u_{1})-f(u_{2})\|~\|A u\|\\
 &\leq L_{f}\|u\|~\|Au\|\\
 &\leq \frac{L_{f}^{2}}{2}\|u\|^2+ \frac{1}{2}\|Au\|^{2}.
 \end{align*}
 Thus,
 \begin{eqnarray}\label{a26}
 \frac{d}{dt}\|u(t)\|^{2} +\|Au(t)\|^{2} \leq c_{f,g}\|u(t)\|^{2},\quad c_{f,g}=(2L_{g}+L_{f}^{2})> 0,\quad t\geq 0.
 \end{eqnarray}
 We finally conclude~\eqref{a27} from~\eqref{a26} and Gronwall's lemma.
 \end{proof}

Next, we show a $L^2-H^1$ smoothing property.
\begin{lem}
\label{a35}
If $\|u_{i}(0)\|_{2}\le  R_{2}$ (i=1,2), then for all $t>0$, we have
\begin{equation}
\label{a35aux}
\|u(t)\|_1^2\le \frac{c_S}{t}\exp(c(R_2)t)\|u_{0}\|^{2}.
\end{equation}	 
\end{lem}
\begin{proof}
We multiply~\eqref{a22} by $2tA^{-1}u_{t}$ in $H$.  We deduce
\begin{eqnarray}
&&\frac{d}{dt}(t\|\nabla u\|^{2}) +2t\|u_{t}\|_{-1}^{2} +2t(g(u_{1})-g(u_{2}),(-\Delta)^{-1}u_{t})+ 2t( f(u_{1})-f(u_{2}),u_{t})\nonumber\\
&&\qquad=\|\nabla u\|^{2}.\label{a36}
\end{eqnarray} 
Using~\eqref{a10}, Poincar\'e's inequality and Young's inequality, we have
\begin{align*}
|(g(u_{1})-g(u_{2}),(-\Delta)^{-1}u_{t})|&\le  \|g(u_{1})-g(u_{2})\|~\|(-\Delta)^{-1}u_{t}\|\\&\le  L_{g}\|u\|~\|u_{t}\|_{-1}\\
&\leqslant L_{g}c_S\|\nabla u\|~\|u_{t}\|_{-1}\\
&\le  L^{2}_{g}c_S^2\|\nabla u\|^{2} + \frac{1}{4}\|u_{t}\|_{-1}^{2}.
\end{align*}
Next, we use the Cauchy-Schwarz inequality and~\eqref{aux2.10}:
\begin{eqnarray}
&&|( f(u_{1})-f(u_{2}),u_{t})|=|(A^{\frac{1}{2}}( f(u_{1})-f(u_{2})),A^{-\frac{1}{2}}u_{t})|\nonumber\\
&&\qquad\leq \left\|\nabla\left(\int_{0}^{1}f'(u_{1}+ s(u_{2}-u_{1}))dsu\right)\right\| \|u_{t}\|_{-1}\nonumber\\
&&\qquad\leq \left\|\int_{0}^{1}f'(u_{1}+ s(u_{2}-u_{1}))ds\nabla u\right\| \|u_{t}\|_{-1}\nonumber\\
&& \qquad\quad+ \left\|\int_{0}^{1}f''(u_{1}+ s(u_{2}-u_{1}))(\nabla u_{1}+ s\nabla(u_{2}-u_{1}))ds u\right\|\|u_{t}\|_{-1}\nonumber\\
&& \qquad\leq  c_f\left(\|\nabla u\|+ \left\| |u|\left|\nabla u_{1}\right|\right\|+ \left\| |u|\left|\nabla u_{2}\right|\right\|\right)\|u_{t}\|_{-1}\nonumber\\
&&\qquad\leq  c_f(\|\nabla u\|+ \|u\|_{L^{4}(\Omega)}\left(\|\nabla u_{1}\|_{L^{4}(\Omega)}+\|\nabla u_{2}\|_{L^{4}(\Omega)}\right)\|u_{t}\|_{-1}\nonumber\\
&&\qquad\leq (\text{thanks ~to~ the~ continuous~ embedding}~H^{1}(\Omega)\subset L^{4}(\Omega))\nonumber\\
&&\qquad\leq  c(\|\nabla u\|+\|\nabla u\|\left(\|u_{1}\|_{H^{2}(\Omega)}+\|u_{1}\|_{H^{2}(\Omega)}\right)\|u_{t}\|_{-1}\nonumber\\
&&\qquad\leq  (\text{since~$u_{1}$,~$u_{2}$ are bounded~in}~H^{2}(\Omega) ~\text{by}~\eqref{a27bis})\nonumber\\
&&\qquad\leq  c(R_2)\|\nabla u\|^{2} +\frac{1}{4}\|u_{t}\|_{-1}^{2}.\label{a45aux}
\end{eqnarray}
Combining the above estimates in~\eqref{a36}, we find
$$\frac{d}{dt}(t\|\nabla u(t)\|^{2}) + t\|u_{t}(t)\|^{2}_{-1}\leq c't\|\nabla u\|^{2} + \| \nabla u(t)\|^{2},\quad c'=2L_{g}^{2}c_S^2+2c(R_2).$$
By Gronwall's lemma, 
$$t\|\nabla u(t)\|^{2}\le e^{c't}\int_0^t\|\nabla u(s)\|^2ds,\quad t\ge 0.$$
We conclude from~\eqref{a27} that~\eqref{a35aux} holds with 
$c(R_2)=c'+c_{f,g}$.
\end{proof}
\section{The time semidiscrete problem}
\label{sec3}
\subsection{The discrete semigroup}
For the time semidiscretization, we apply the semi-implicit Euler scheme to~\eqref{a1}.
In the remainder of the manuscript, $\delta t >0$ denotes the time step. The scheme reads : let $u^{0}\in D(A)=H^2(\Omega)\cap H^1_0(\Omega)$ and for $n=0,1,2,\cdots, $ let $u^{n+1} \in D(A)$ solve
\begin{equation}
\label{b1}
\frac{u^{n+1}-u^{n}}{\delta t} +A^{2} u^{n+1} +A f(u^{n})+ g(u^{n})=0.
\end{equation}
This is a linear sheme known as the IMEX (Implicit-Explicit) scheme: at each time step, $u^{n+1}$ is computed by solving a linear system whose right-hand side involves $u^n$.   By elliptic regularity, for $u^{n}\in D(A)$, $u^{n+1}$ is unique and belongs to $D(A)$.
The following result shows that the discrete semigroup $S^{n}_{\delta t}u^{0}=u^{n}$ is well-defined on $D(A)$. 
\begin{theo}\label{b2}
Assume that $\delta t\leq 1/(2L_{g})$, where $L_g$ is the constant in~\eqref{c12bis}. Then for every $u^{n}\in D(A)$, there exists a unique $u^{n+1}\in D(A)$ which solves~\eqref{b1}. 
Moreover, the mapping $S_{\delta t}:u^{n}\mapsto u^{n+1}$ is  Lipschitz continuous for the $L^2(\Omega)$-norm from $D(A)$ into $D(A)$.
\end{theo}
The Lipschitz continuity in $L^2(\Omega)$ follows from~\eqref{c2aux}. 

The following regularity result will prove useful:
\begin{lem}
\label{lemb2}
If $u^n\in D(A)$, then $u^{n+1}=S_{\delta t}u^n$ belongs to $D(A^2)$ and 
\begin{equation}
\label{eqA2dt}
\delta t\|A^2u^{n+1}\|^2\le 2\|\Delta u^n\|^2+C\delta t\left(\|\Delta u^n\|^4+1\right),
\end{equation}
where the positive constant $C$ is independent of $\delta t$ and $u^n$. 
Moreover, 
\begin{equation}
\label{eqDA}
\|\Delta u^{n+1}\|^2\le \|\Delta u^n\|^2+C\delta t\left(\|\Delta u^n\|^4+1\right).
\end{equation}
\end{lem}
\begin{proof}By~\eqref{b1}, $u^{n+1}$ solves
\begin{equation}
\label{aux1}
u^{n+1}-u^n+\delta t A^2u^{n+1}=\delta t\, h
\end{equation}
where $h=\Delta f(u^n)-g(u^n)$. By the chain rule~\eqref{z1}, $h\in L^2(\Omega)$ with
\begin{eqnarray}
\|h\|&\le &\|f^\prime(u^n)\|_{L^\infty(\Omega)}\|\Delta u^n\|+\|f^{\prime\prime}(u^n)\|_{L^\infty(\Omega)}\|\nabla u^n\|^2_{L^4(\Omega)}+\|g(u^n)\|\nonumber\\
&\le & (H^1(\Omega)\subset L^4(\Omega))\nonumber\\
&\le & c_{f'}\|\Delta u^n\|+c_{f''}c_S^2\|\Delta u^n\|^2+c_g .\label{aux3}
\end{eqnarray}
Since $u^n$, $u^{n+1}$ and $h$ belong to $L^2(\Omega)$, we deduce from~\eqref{aux1} that $u^{n+1}\in D(A^2)$. Next, we take the $L^2$-scalar product of~\eqref{aux1} with $u^{n+1}-u^n$. This yields
\begin{eqnarray*}\|u^{n+1}-u^n\|^2+\delta t\|\Delta u^{n+1}\|^2&=&\delta t(\Delta u^{n+1},\Delta u^n)+\delta t(h,u^{n+1}-u^n)\\
&\le &\delta t\|\Delta u^{n+1}\|\,\|\Delta u^n\|+\delta t\|h\|\,\|u^{n+1}-u^n\|.
\end{eqnarray*}
From Young's inequality, we deduce that 
\begin{equation}
\label{aux2}
\|u^{n+1}-u^n\|^2+\delta t\|\Delta u^{n+1}\|^2\le \delta t\|\Delta u^n\|^2+\delta t^2\|h\|^2.
\end{equation}
Now, we take the $L^2$-scalar product of~\eqref{aux1} with $A^2u^{n+1}$. We find
\begin{eqnarray*}
\delta t\|A^2u^{n+1}\|^2&=&-(u^{n+1}-u^n,A^2u^{n+1})+\delta t(h,A^2u^{n+1})\\
&\le &\frac{1}{\delta t}\|u^{n+1}-u^n\|^2+\frac{\delta t}{4}\|A^2u^{n+1}\|^2+\delta t\|h\|^2+\frac{\delta t}{4}\|A^2u^{n+1}\|^2.
\end{eqnarray*}
Thus, by~\eqref{aux2}, 
$$\delta t\|A^2u^{n+1}\|^2\le 2\|\Delta u^n\|^2+4\delta t\|h\|^2.$$
Using~\eqref{aux3}, we find~\eqref{eqA2dt}. From~\eqref{aux2} and~\eqref{aux3}, we also deduce that
\begin{eqnarray*}
\|\Delta u^{n+1}\|^2&\le &\|\Delta u^n\|^2+\delta t\|h\|^2\\
&\le &\|\Delta u^n\|^2+C\delta t\left(\|\Delta u^n\|^4+1\right).
\end{eqnarray*}
This is~\eqref{eqDA}. 
\end{proof}
\subsection{Dissipative estimates, uniform in $\delta t$}
\label{subsec3.2}
In this subsection, we first establish some priori estimates for the scheme~\eqref{b1}. The following well-known identity will be frequently used:
\begin{equation}\label{b5}
(a-b,a)=\frac{1}{2}(\|a\|^{2}- \|b\|^{2}+ \|a-b\|^{2}),\quad a,b \in L^{2}(\Omega)
\end{equation}
In~\eqref{b5}, we may also replace the  inner product and the norm in  $L^2(\Omega)$ by another inner product and the norm associated to it. 
We  recall a discrete Gronwall lemma.
\begin{lem}
\label{lemGD}
Let $C,\gamma>0$ and $(a_n)$, $(b_n)$ be two sequences of nonnegative real numbers such that 
\begin{equation}
\label{auxb1}
a_{n+1}+\delta t b_{n+1}\le (1-\gamma\delta t)a_n+\delta t C,\quad\forall n\ge 0,
\end{equation}
where $\delta t\in(0,1/(2\gamma)]$. 
Then for all $n\ge 0$, we have
\begin{equation}
\label{auxb2}
a_n\le e^{-n\gamma\delta t}a_0+\frac{C}{\gamma}
\end{equation}
and 
\begin{equation}
\label{auxb3}
\delta t\sum_{k=0}^{n-1}b_{k+1}\le a_0+n\delta t C.
\end{equation}
\end{lem}
By convention, for $n=0$, the sum in the left-hand side of~\eqref{auxb3} is zero. 
\begin{proof}Since
$$a_{n+1}\le (1-\gamma\delta t)a_n+\delta t C,\quad\forall n\ge 0,$$
we find by induction that for all $n\ge 0$, 
\begin{eqnarray*}
a_n&\le& (1-\gamma\delta t)^n a_0+\delta tC\sum_{k=0}^{n-1}(1-\gamma\delta t)^k\\
&\le & (1-\gamma\delta t)^n a_0+\frac{C}{\gamma}.
\end{eqnarray*}
By convexity, we have $1-s\le e^{-s}$ for all $s\ge 0$. Thus, for all $n\ge 0$, 
$$a_n\le e^{-n\gamma\delta t}a_0+\frac{C}{\gamma}.$$
This is~\eqref{auxb2}.  By~\eqref{auxb1}, we also have 
$$a_{n+1}+\delta t b_{n+1}\le a_n+\delta tC,\quad\forall n\ge 0.$$
By summing from $n=0$ to $n=N-1$, we find 
$$a_N+\delta t\sum_{n=0}^{N-1}b_{n+1}\le a_0+N\delta t C.$$
This yields~\eqref{auxb3}. 
\end{proof}
\begin{prop}[Dissipative estimate  in $L^{2}(\Omega)$]
\label{b6}
If $\delta t$ is small enough, then
\begin{equation}
\label{b6bis}
\|u_n\|^2\le C_0\|u_0\|^2e^{-\epsilon_0 n \delta t}+M_0,\quad \forall n\ge 0,
\end{equation}
and 
\begin{equation}
\label{b6ter}
\delta t\sum_{k=0}^{n-1}\|\Delta u_{k+1}\|^2\le C_0'\|u_0\|^2+n\delta tM_0',\quad \forall n\ge 0,
\end{equation}
where the positive constants  $C_0$, $\epsilon_0$, $M_0$, $C_0'$  and $M_0'$ are independent of $u_0$ and $\delta t$. 
\end{prop}
\begin{proof}
We multiply~\eqref{b1} by $A^{-1}u^{n+1}$ in $H$. Using~\eqref{b5}, we have
\begin{eqnarray}\label{b8}
\frac{1}{2\delta t}(\|u^{n+1}\|^{2}_{-1}-\|u^{n}\|^{2}_{-1}+\|u^{n+1}-u^{n}\|^{2}_{-1})+ \|\nabla u^{n+1}\|^{2} +(f(u^{n+1}),u^{n+1})\nonumber\\+(f(u^{n})-f(u^{n+1}),u^{n+1})+ (g(u^{n}),(-\Delta)^{-1}u^{n+1})=0. 
\end{eqnarray}
Thanks to~\eqref{a7bis}, we have
\begin{align*}
(f(u^{n+1}),u^{n+1})\geq c_{1}  \|u^{n+1}\|^{2}-c_{2},\quad c_{1}=\gamma_5>0,~c_{2}\geq 0.
\end{align*}
Using the Cauchy-Schwarz inequality, ~\eqref{c12f} and Young's inequality, we find
\begin{align*}
|(f(u^{n})-f(u^{n+1}),u^{n+1})|&\leq \|f(u^{n})-f(u^{n+1})\|~\|u^{n+1}\|\\&\leq L_{f}\|u^{n}-u^{n+1}\|~\|u^{n+1}\|\\&\leq \frac{c_{1}}{4}\|u^{n+1}\|^{2}+ \frac{L^{2}_{f}}{c_{1}}\|u^{n}-u^{n+1}\|^{2}.
\end{align*}
Moreover, 
\begin{align*}
|(g(u^{n}),(-\Delta)^{-1}u^{n+1})|&\leq \|g(u^{n})\|~\|(-\Delta)^{-1}u^{n+1}\|\\&\leq c_{g} c_{S}\|u^{n+1}\|\\
&\leq \frac{c_{1}}{4}\|u^{n+1}\|^{2}+ \frac{c'^{2}_{g}}{c_{1}},\quad \mbox{ where }c_g'=c_gc_S. 
\end{align*}
Let us combine the above estimates in~\eqref{b8}. We find
\begin{eqnarray}
&&\frac{1}{2\delta t}\|u^{n+1}\|_{-1}^{2} +\frac{c_{1}}{2}\|u^{n+1}\|^{2} + \|\nabla u^{n+1}\|^{2}+ \frac{1}{2\delta t}\|u^{n+1}-u^{n}\|^{2}_{-1}\nonumber\\
&&\leq \frac{1}{2\delta t}\|u^{n}\|^{2}_{-1}+\frac{L^{2}_{f}}{c_{1}}\|u^{n+1}-u^{n}\|^{2}+\frac{c'^{2}_{g}}{c_{1}}.\label{b9}
\end{eqnarray}

Now, we multiply~\eqref{b1} by $u^{n+1}$ in $H$.   We obtain
\begin{eqnarray}\label{b10}
&&\frac{1}{2\delta t}(\|u^{n+1}\|^{2}-\|u^{n}\|^{2}+\|u^{n+1}-u^{n}\|^{2})+ \|\Delta u^{n+1}\|^{2}\nonumber\\
&&+(f(u^{n}),\Delta u^{n+1})+ (g(u^{n}),u^{n+1}) =0.
\end{eqnarray}
By~\eqref{c12f},~\eqref{a4},~\eqref{sobo} and Young's inequality, we have
\begin{align*}
|(f(u^{n}),\Delta u^{n+1})|&\leq \|f(u^{n})\|~\|\Delta u^{n+1}\|\\
&\leq L_{f}\|u^{n}\|~\|\Delta u^{n+1}\|\\
&\leq \frac{L^{2}_{f}c_S^2}{2}\|\nabla u^{n}\|^{2}+ \frac{1}{2}\|\Delta u^{n+1}\|^{2}.
\end{align*}
Owing to~\eqref{c12} and  Young's inequality, we have
\begin{align*}
|(g(u^{n}),u^{n+1})|&\leq \|g(u^{n})\|~\|u^{n+1}\|\\
&\leq c_{g}\|u^{n+1}\|\\
&\leq \frac{1}{4}\|u^{n+1}\|^{2}+c^{2}_{g}.
\end{align*}
Combining the above estimates in~\eqref{b10}, we obtain
\begin{eqnarray}
&&\frac{1}{2\delta t}\|u^{n+1}\|^{2}+\frac{1}{2 \delta t}\|u^{n+1}-u^{n}\|^{2}+\frac{1}{2}\|\Delta u^{n+1}\|^{2}\nonumber\\
&&\qquad\leq \frac{1}{2 \delta t}\|u^{n}\|^{2}+\frac{L^{2}_{f}c_S^2}{2}\|\nabla u^{n}\|^{2}+ \frac{1}{4}\|u^{n+1}\|^{2}+ c^{2}_{g}.\label{b11}
\end{eqnarray}
On summing~\eqref{b9} and $\alpha$~\eqref{b11} with $\alpha>0$ small enough, we conclude that
\begin{align}
&\frac{1}{2\delta t}\|u^{n+1}\|_{-1}^{2}+\frac{\alpha}{2\delta t}\|u^{n+1}\|^{2} +(\frac{c_{1}}{2}-\frac{\alpha}{4})\|u^{n+1}\|^{2} + \|\nabla u^{n+1}\|^{2}\nonumber\\
&+ \frac{1}{2\delta t}\|u^{n+1}-u^{n}\|^{2}_{-1}+ \frac{\alpha}{2 \delta t}\|u^{n+1}-u^{n}\|^{2}+\frac{\alpha}{2}\|\Delta u^{n+1}\|^{2}\nonumber\\
&\qquad \leq \frac{1}{2\delta t}\|u^{n}\|^{2}_{-1} +\frac{\alpha}{2 \delta t}\|u^{n}\|^{2}+\frac{L^{2}_{f}}{c_{1}}\|u^{n+1}-u^{n}\|^{2} +\alpha \frac{L^{2}_{f}c_S^2}{2}\|\nabla u^{n}\|^{2} + c', \label{b12}
\end{align}
where $c'=\cfrac{c_g^{\prime 2}}{c_{1}}+\alpha c_g^2$.
We choose $\alpha>0$ small enough so that $\alpha\le c_{1}$ and $\alpha L_f^2c_S^2/2\le 1$ and  
we set 
\begin{equation}
\label{defan}
a_{n}=\|u^{n}\|^{2}_{-1}+\alpha\|u^{n}\|^{2}.
\end{equation}
 Then,  for $\delta t$ small enough (depending only on $\alpha$ and $L_f^2/c_{1}$),  the estimate~\eqref{b12} yields
\begin{eqnarray*}
\frac{1}{2\delta t}a_{n+1}+ \frac{c_{1}'}{2}a_{n+1}+\frac{\alpha}{2}\|\Delta u^{n+1}\|^2\leq \frac{1}{2\delta t}a_{n}+ c',\quad c_{1}'=c_{1}'(c_{1},c_S)>0.
\end{eqnarray*}
Thus,
\begin{eqnarray*}
a_{n+1}+\frac{\alpha \delta t}{1+c_{1}'\delta t}\|\Delta u^{n+1}\|^2\leq \frac{1}{1+c_{1}'\delta t}a_{n}+ \frac{2\delta t}{1+c_{1}'\delta t} c'.
\end{eqnarray*}
Thanks to an asymptotic expansion of   order 2, for $\delta t$ small enough (depending only on $c_{1}'$), we have
\begin{equation}
\label{b14}
a_{n+1}+\frac{\alpha}{2}\delta t\|\Delta u^{n+1}\|^2\leq (1-\frac{c_{1}'}{2}\delta t)a_{n}+ 2c'\delta t.
\end{equation}
We deduce from Lemma~\ref{lemGD} that for all $n\ge 0$, 
$$a_n\le e^{-c_{1}'n\delta t/2}a_0+\frac{4c'}{c_{1}'}$$
 and 
 $$\frac{\alpha}{2} \delta t\sum_{k=0}^{n-1}\|\Delta u^{k+1}\|^2\le a_0+2c'n\delta t.$$
 This concludes the proof of Proposition~\ref{b6}.
\end{proof}

\begin{prop}[Dissipative estimate in $H^{1}(\Omega)$] 
If $\delta t$ is small enough, then
\begin{equation}
\label{auxb7}
\|\nabla u^n\|^2\le C_{1}e^{-\epsilon_1n\delta t}\|\nabla u^0\|^2+M_1,\quad\forall n\ge 0,
\end{equation}
and 
$$\delta t\sum_{k=0}^{n-1}\|u^{k+1}-u^k\|_{-1}^2\le C_{1}'\|\nabla u^0\|^2+M_1', \quad\forall n\ge 0,$$
where the positive constants $C_{1}$, $\epsilon_1$, $M_1$, $C_{1}'$ and $M_1'$ are independent of $u_0$ and $\delta t$. 
\end{prop}
\begin{proof}
We multiply~\eqref{b1} by $A^{-1}\cfrac{u^{n+1}-u^{n}}{\delta t}$ in $H$.  We obtain 
\begin{eqnarray}\label{b17}
&&\frac{1}{\delta t ^{2}}\|u^{n+1}-u^{n}\|^{2}_{-1}+\frac{1}{2\delta t}\left(\|\nabla u^{n+1}\|^{2}-\|\nabla u^{n}\|^{2}+\|\nabla(u^{n+1}-u^{n})\|^{2}\right)\nonumber\\
&&+\frac{1}{\delta t}(F(u^{n+1})-F(u^{n}),1)+ ( g(u^{n}) , (-\Delta)^{-1}\frac{u^{n+1}-u^{n}}{\delta t})\nonumber\\
&&\qquad= \frac{1}{2\delta t}\int_{\Omega}f'(\zeta_{u^{n+1},u^{n}})(u^{n+1}-u^{n})^{2}dx.
\end{eqnarray}
Here, we used that for all $r,s\in\Rr$, 
$$F(s)=F(r)+f(r)(s-r)+f'(\xi_{s,r})\frac{(s-r)^2}{2},\ \mbox{ for some } \xi_{s,r}\in[r,s],$$
and $(r,s)\mapsto f'(\xi_{s,r})$ is a continuous function on $\Rr^2$, since $F\in C^2(\Rr)$.  

Using~\eqref{a5}, an interpolation inequality and Young's inequality, we obtain
\begin{align*}
\left|\int_{\Omega}f'(\zeta_{u^{n+1},u^{n}})(u^{n+1}-u^{n})^{2}dx\right|&\leq
\|f'\|_{L^\infty(\Omega)}~\|u^{n+1}-u^{n}\|^{2}\\
&\leq L_{f}\|u^{n+1}-u^{n}\|^{2}\\
&\leq L_f\|\nabla (u^{n+1}-u^n)\|\,\|u^{n+1}-u^n\|_{-1}\\
&\leq\frac{1}{2}\|\nabla(u^{n+1}-u^{n})\|^{2}+\frac{L^{2}_{f}}{2}\|u^{n+1}-u^{n}\|^{2}_{-1}.
\end{align*}
By~\eqref{c12},
\begin{align*}
\left|( g(u^{n}) , (-\Delta)^{-1}\frac{u^{n+1}-u^{n}}{\delta t})\right|&\leq \| g(u^{n})\|\left\| (-\Delta)^{-1}\frac{u^{n+1}-u^{n}}{\delta t}\right\|\\
&\leq c_{g}c_S\|\frac{u^{n+1}-u^{n}}{\delta t}\|_{-1}\\
&\leq \frac{1}{2\delta t^{2}}\|u^{n+1}-u^{n}\|_{-1}^{2}+ \frac{c^{2}_{g}c_S^2}{2}.
\end{align*}
Combining the above estimates in~\eqref{b17}, we get
\begin{eqnarray}
&&\frac{1}{2\delta t}\|\nabla u^{n+1}\|^{2}+ \frac{1}{\delta t}(F(u^{n+1}),1)+ \frac{1}{2\delta t^{2}}\|u^{n+1}-u^{n}\|_{-1}^{2}+\frac{1}{4\delta t}\|\nabla(u^{n+1}-u^{n})\|^{2}\nonumber \\
&&\qquad\leq \frac{1}{2\delta t}\|\nabla u^{n}\|^{2}+\frac{1}{\delta t}(F(u^{n}),1)+\frac{L^{2}_{f}}{4\delta t}\|u^{n+1}-u^{n}\|_{-1}^{2}+ \frac{c^{2}_{g}c_S^2}{2}.\label{b18}
\end{eqnarray}
Let $\delta t$ be small enough so that~\eqref{b14} holds and $\delta t\le 1/(4L_f^2)$. 
Adding~\eqref{b14}  and $2\delta t\beta$ times~\eqref{b18}  where $\beta>0$ is small enough, we find that
\begin{equation}
\label{b20}
E_2^{n+1}+\frac{\alpha}{2}\delta t\|\Delta u^{n+1}\|^2+\frac{\beta}{2\delta t}\|u^{n+1}-u^n\|_{-1}^2\le E_2^n+(2c'+\beta c_g^2c_S^2)\delta t,
\end{equation}
where 
$$E_2^n=a_n+\beta\|\nabla u^n\|^2+2\beta (F(u^n),1)$$
and $a_n$ is defined by~\eqref{defan}. By~\eqref{e1}, $F(s)+\gamma_2\ge 0$ for all $s\in\Rr$, so that 
$$\tilde{E}_2^n=E_2^n+2\beta(\gamma_2,1)\ge 0,\quad\forall n\ge 0.$$
Moreover, by~\eqref{e2} and the Poincar\'e inequality, we have
\begin{equation}
\label{auxb3bis}
\tilde{E}_2^n\le c_{3}\|\nabla u^n\|^2+c_{4},\quad c_{3},\  c_{4}>0,
\end{equation}
so that~\eqref{b20} yields
$$\tilde{E}_2^{n+1}+c_{5}\delta t\tilde{E}_2^{n+1}+\frac{\beta}{2\delta t}\|u^{n+1}-u^n\|_{-1}^2\le \tilde{E}_2^n+c_{6}\delta t,\quad c_{5},\ c_{6}>0.$$
Thus, 
$$\tilde{E}_2^{n+1}+\frac{1}{1+c_{5}\delta t}\frac{\beta}{2\delta t}\|u^{n+1}-u^n\|_{-1}^2\le\frac{1}{1+c_{5}\delta t} \tilde{E}_2^n+\frac{c_{6}}{1+c_{5}\delta t}\delta t.$$
Therefore, for $\delta t$ small enough (depending only on $c_{5}$), we have 
$$\tilde{E}_2^{n+1}+\frac{\beta\delta t}{4}\left\|\frac{u^{n+1}-u^n}{\delta t}\right\|_{-1}^2\le (1-\frac{c_{5}}{2}\delta t) \tilde{E}_2^n+c_{6}\delta t.$$
We may apply Lemma~\ref{lemGD}, which yields
\begin{equation}
\label{auxb4}
\tilde{E}_2^n\le e^{-nc_{5}\delta t/2}\tilde{E}_2^0+\frac{2c_{6}}{c_{5}}
\end{equation}
and 
\begin{equation}
\label{auxb5}
\frac{\beta}{4\delta t}\sum_{k=0}^{n-1}\|u^{k+1}-u^k\|_{-1}^2\le \tilde{E}_2^0+n\delta t c_{6},
\end{equation}
for all $n\ge 0$. 
Finally, we note that 
$$\tilde{E}_2^n\ge \beta\|\nabla u^n\|^2,$$
and this estimate, together with~\eqref{auxb3bis}, \eqref{auxb4} and~\eqref{auxb5}, concludes the proof. 
\end{proof}

\begin{prop}[Dissipative estimate in $H^2(\Omega)$]
\label{propH2dt}
 For $\delta t$ small enough, we have
\begin{equation}
\label{dissiH2dt}
\|\Delta u^n\|^2\le Q_2(\|\Delta u^0\|)e^{-\epsilon_2 n\delta t}+M_2,\quad \forall n\ge 0,
\end{equation}
and 
\begin{equation}
\label{dissiH2dtbis}
\sum_{k=0}^{n-1}\|\Delta (u^{k+1}-u^k)\|^2\le Q_2(\|\Delta u^0\|)+M_2'n\delta t,\quad n\ge 0,
\end{equation}
where the monotonic function $Q_2$ and the positive constants $\epsilon_2$, $M_2$ and $M_2'$ are independent of $u^0$ and $\delta t$.
\end{prop}
\begin{proof}
By Lemma~\ref{lemb2}, we know that for all $n\ge 1$, $u^{n}\in D(A^2)$. 
On multiplying~\eqref{b1} by $A^{2}u^{n+1}$ in $H$ and using~\eqref{b5}, we obtain
\begin{eqnarray}
&&\frac{1}{2\delta t}\left(\|\Delta u^{n+1}\|^{2}-\|\Delta u^{n}\|^{2}+\|\Delta(u^{n+1}-u^{n})\|^{2}\right)+\|\Delta^{2}u^{n+1}\|^{2}\nonumber\\
&&=-(g(u^{n}),\Delta^{2}u^{n+1})+(\Delta f(u^{n}),\Delta^{2}u^{n+1}).\label{b24}
\end{eqnarray} 
Using~\eqref{c12} and Young's inequality yields
\begin{align*}
|(g(u^{n}),\Delta^{2}u^{n+1})|&\leq\|g(u^{n})\|~\|\Delta^{2}u^{n+1}\|\\&
\leq c_{g}\|\Delta^{2}u^{n+1}\|\\&
\leq \frac{1}{4}\|\Delta^{2}u^{n+1}\|^{2}+c^{2}_{g}.
\end{align*}
Similarly, we have
\begin{align*}
|(\Delta f(u^{n}),\Delta^{2}u^{n+1})|&\leq \|\Delta f(u^{n})\|~\|\Delta^{2}u^{n+1}\|\\
&\leq \frac{1}{4}\|\Delta^{2}u^{n+1}\|^{2}+\|\Delta f(u^{n})\|^{2}.
\end{align*}
Therefore, we have
\begin{eqnarray}
&&\frac{1}{2\delta t}\left(\|\Delta u^{n+1}\|^{2}-\|\Delta u^{n}\|^{2}+\|\Delta(u^{n+1}-u^{n})\|^{2}\right)+\frac{1}{2}\|\Delta^{2}u^{n+1}\|^{2}\nonumber\\
&&\qquad\qquad\leq\|\Delta f(u^{n})\|^{2}+c^{2}_{g}.\label{b25}
\end{eqnarray} 
Arguing as in the continuous case (cf.~\eqref{a30bis}) and using the $H^1$-estimate~\eqref{auxb7}, we find that for  all $n\ge 1$, we have
\begin{equation}
\label{auxc19}
\|\Delta f(u^n)\|\le \frac{1}{4}\|\Delta^2u^n\|^2+c_{16}(R_1),
\end{equation}
where $c_{16}$ depends on $R_1=\|\nabla u(0)\|$. 
The estimate~\eqref{b25}  becomes
\begin{eqnarray}
&&\frac{1}{2\delta t}\|\Delta u^{n+1}\|^{2}+\frac{1}{2}\|\Delta^{2}u^{n+1}\|^{2}+\frac{1}{2\delta t}\|\Delta(u^{n+1}-u^{n})\|^{2}\nonumber\\
&&\qquad\leq\frac{1}{2\delta t}\|\Delta u^{n}\|^{2} +\frac{1}{4}\|\Delta^{2}u^{n}\|^{2}+c_{7},\nonumber
\end{eqnarray}
where $c_{7}(R_1)=c_g^2+c_{16}(R_1)$. We multiply this estimate by $2\delta t$ and we obtain
\begin{equation}
\label{aux4}
b_{n+1}+\frac{\delta t}{2}\|\Delta^2u^{n+1}\|^2+\|\Delta(u^{n+1}-u^n)\|^2\le b_n+2 c_{7}\delta t,
\end{equation}
where 
$$b_n=\|\Delta u^n\|^2+\frac{\delta t}{4}\|\Delta^2 u^{n}\|^2.$$
Using that  $\delta t$ is bounded from above and that~\eqref{sobo} holds for $s_1=1$ and $s_2=2$, we find
$$b_n\le c_{8}\|\Delta^2u^n\|^2,\quad\forall n\ge 1$$
for some constant $ c_{8}>0$ independent of $\delta t$. 
Thus,~\eqref{aux4} implies
$$b_{n+1}+\frac{\delta t}{2c_{8}}b_{n+1}+\|\Delta(u^{n+1}-u^n)\|^2\le b_n+2c_{7}\delta t,\quad\forall n\ge 1.$$
Therefore, for $\delta t$ small enough, we have
\begin{equation}
\label{auxbn}
b_{n+1}+\frac{1}{2}\|\Delta(u^{n+1}-u^n)\|^2\le (1-c_{9}\delta t)b_n+2c_{7}\delta t\quad\forall n\ge 1,
\end{equation}
where $c_{9}=1/(4c_{8})$. By Lemma~\ref{lemGD}, 
\begin{equation}
\label{auxbnter}
b_n\le e^{-(n-1)c_{9}\delta t}b_1+\frac{2c_{7}}{c_{9}},\quad\forall n\ge 1,
\end{equation}
and
\begin{equation}
\label{auxbn5}
\frac{1}{2}\sum_{k=1}^{n-1}\|\Delta (u^{k+1}-u^k)\|^2\le b_1+(n-1)\delta t\frac{2c_{7}}{c_{9}},\quad\forall n\ge 1.
\end{equation}

Thanks to the dissipative estimate~\eqref{auxb7}, there exists a time $t_1=t_1(R_1)$ such that 
$$\|\nabla u^n\|^2\le C_{1}+M_1,\quad\forall n\ge t_1/\delta t.$$
Let $n_1=\lceil t_1/\delta t \rceil$, where $\lceil\cdot\rceil$ denotes the integer ceiling function. Then for $n\ge n_1$, the constant $c_{16}$ in~\eqref{auxc19} no longer depends on $R_1$. Consequently,~\eqref{auxbn} holds for all $n\ge  n_1$ with $c_{7}$ and $c_{9}$ independent of $R_1$. By induction (as in Lemma~\ref{lemGD}), we have
\begin{equation}
\label{auxbnbis}
b_n\le e^{-(n-n_1)c_{9}\delta t} b_{n_1}+\frac{2c_{7}}{c_{9}},\quad\forall n\ge n_1,
\end{equation}
and
\begin{equation}
\label{auxbn4}
\frac{1}{2}\sum_{k=n_1}^{n-1}\|\Delta (u^{k+1}-u^k)\|^2\le b_{n_1}+(n-n_1)\delta t\frac{2c_{7}}{c_{9}},\quad\forall n\ge n_1.
\end{equation}
The dissipative estimate~\eqref{dissiH2dt} follows from~\eqref{auxbnbis} (for $n\ge n_1$), \eqref{auxbnter} (for $1\le n\le n_1$ and Lemma~\ref{lemb2} (for $n=0$). Estimate~\eqref{dissiH2dtbis} follows from~\eqref{auxbn4} (for $n\ge n_1$), \eqref{auxbn5} (for $1\le n\le n_1$), \eqref{auxbnter} (for $b_{n_1}$) and Lemma~\ref{lemb2} (for $n=0$). 
 \end{proof}		
\subsection{Estimates for the difference of solutions, uniform in $\delta t$}
Let $v^{n}$ and $w^{n}$ be two sequences generated by the scheme~\eqref{b1}  and corresponding to the initial data $v^{0}$ and $w^{0}$ respectively. We denote $u^{n}=v^{n}-w^{n}$ their difference, which satisfies
\begin{eqnarray}\label{c1}
\frac{u^{n+1}-u^{n}}{\delta t} +A^{2}u^{n+1}+A (f(v^{n})-f(w^{n}))+( g(v^{n})-g(w^{n}))=0,\quad \forall n\geq 0.
\end{eqnarray}
\begin{lem}\label{c2}
Assume that $\delta t<1/(2L_g)$. Then for all $n\ge 1$, we have
\begin{equation}
\label{c2aux}
\|u^{n}\|^{2}+ \sum_{k=0}^{n-1}\|u^{k+1}-u^{k}\|^{2}+\delta t \sum_{k=0}^{n-1}\|\Delta u^{k+1}\|^{2}\leq \exp{(c_{f,g}n\delta t)}\|u^{0}\|^{2}.
\end{equation}
\end{lem}
\begin{proof}
We multiply~\eqref{c1}  by $u^{n+1}$ in $H$. We obtain
\begin{eqnarray}
&&\frac{1}{2\delta t}(\|u^{n+1}\|^{2}-\|u^{n}\|^{2}+\|u^{n+1}-u^{n}\|^{2})+\|\Delta u^{n+1}\|^{2}\nonumber\\
&&\qquad\qquad=(f(v^{n})-f(w^{n}),\Delta u^{n+1})- (g(v^{n})-g(w^{n}),u^{n+1}).\label{c3}
\end{eqnarray}
Owing  to~\eqref{c12bis} and  Young's inequality, we have
\begin{align*}
|(g(v^{n})-g(w^{n}),u^{n+1})|&\leq \|g(v^{n})-g(w^{n})\|~\|u^{n+1}\|\\&
\leq L_{g}\|u^{n}\|~\|u^{n+1}\|\\
&\leq \frac{L_{g}}{2}\|u^{n+1}\|^{2}+\frac{L_g}{2}\|u^{n}\|^{2}.
\end{align*}
By~\eqref{c12f} and Young's inequality,
\begin{align*}
|(f(v^{n})-f(w^{n}),\Delta u^{n+1})|&\leq \|f(v^{n})-f(w^{n})\|~\|\Delta u^{n+1}\|\\
&\leq L_{f}\|u^{n}\|~\|\Delta u^{n+1}\|\\
&\leq \frac{L^{2}_{f}}{2}\|u^{n}\|^{2}+ \frac{1}{2}\|\Delta u^{n+1}\|^{2}.
\end{align*}
Plugging this in~\eqref{c3} times $2\delta t$, we find
\begin{equation}\label{c4}
(1-L_g\delta t)\|u^{n+1}\|^{2}+ \|u^{n+1}-u^{n}\|^{2}+\delta t\|\Delta u^{n+1}\|^{2}\leq (1+c\delta t)\|u^{n}\|^{2},
\end{equation}
where $c=L_f^2+L_g$. 
Thus, for $\delta t\le 1/(2L_g)$, we have
\begin{eqnarray}\label{c5}
\|u^{n+1}\|^{2}+ \|u^{n+1}-u^{n}\|^{2}+\delta t\|\Delta u^{n+1}\|^{2}\leq (1+c'\delta t)\|u^{n}\|^{2},\quad \forall n\geq 0,
\end{eqnarray}
where $c'=c'(c,L_g)$. We apply the estimate
\begin{equation}\label{c21}
1+s\leq \exp{(s)},\quad\forall s\in\mathbb{R},
\end{equation}
to $s=c'\delta t$  and we obtain~\eqref{c2aux}  by induction, with $c_{f,g}=c'$. 
\end{proof}

Next, we show a $L^2$-$H^1$ smoothing property.
\begin{lem}
\label{lem30}
Let $R_{2}>0$ and $\delta t<1/(2L_g)$. If $\|v^{0}\|_2\leq R_{2}$ and $\| w^{0}\|_2\leq R_{2}$, then for all $n\geq 1$, we have
\begin{eqnarray}\label{c6}
n\delta t\|u^{n}\|_{1}^{2}\leq c_S\exp(c(R_2)n\delta t)\|u^{0}\|^{2}.
\end{eqnarray}
\end{lem}
\begin{proof}
 We multiply~\eqref{c1}  by $A^{-1}(u^{n+1}-u^{n})/\delta t$ in $H$ and we find 
 \begin{eqnarray}
&& \frac{1}{\delta t^{2}}\|u^{n+1}-u^{n}\|^{2}_{-1}+ \frac{1}{2\delta t}\left(\|\nabla u^{n+1}\|^{2}-\|\nabla u^{n}\|^{2}+\|\nabla(u^{n+1}-u^{n}\right)\|^{2})\nonumber\\
&&+\left(f(v^{n})-f(w^{n}),\frac{u^{n+1}-u^{n}}{\delta t}\right)\nonumber\\
&&+ \left(g(v^{n})-g(w^{n}),A^{-1}\frac{u^{n+1}-u^{n}}{\delta t}\right)=0.\label{c7}
 \end{eqnarray}
 Using~\eqref{c12bis}, \eqref{sobo}, the Poincar\'e inequality and Young's inequality, we get
 \begin{align*}
 |(g(v^{n})-g(w^{n}),(-\Delta)^{-1}\frac{u^{n+1}-u^{n}}{\delta t}) |&\leq \|g(v^{n})-g(w^{n})\|~\|(-\Delta)^{-1}\frac{u^{n+1}-u^{n}}{\delta t}\|\\
 &\leq L_{g}\|u^{n}\|\,c_S\|\frac{u^{n+1}-u^{n}}{\delta t}\|_{-1}\\
 &\leq L_g^2c_S^2\|\nabla u^{n}\|^{2}+\frac{1}{4\delta t^{2}}\|u^{n+1}-u^{n}\|_{-1}^{2}.
 \end{align*}
 Thanks to~\eqref{dissiH2dt}, we know that $(v^n)$ and $(w^n)$ are bounded in $H^2(\Omega)$. 
 Arguing as in the continuous case (see~\eqref{a45aux}), we obtain
 \begin{align*}
 |(f(v^{n})-f(w^{n}),\frac{u^{n+1}-u^{n}}{\delta t}) |&\leq \|A^{\frac{1}{2}}(f(v^{n})-f(w^{n}))\|~\|\frac{u^{n+1}-u^{n}}{\delta t}\|_{-1}\\
 &\leq c(R_2)\|\nabla u^{n}\|^{2}+\frac{1}{4\delta t^{2}}\|u^{n+1}-u^{n}\|_{-1}^{2}.
 \end{align*}
 We combine the above estimates in~\eqref{c7} and we deduce that
$$ \frac{1}{2\delta t}\|\nabla u^{n+1}\|^{2}\leq \frac{1}{2\delta t}\|\nabla u^{n}\|^{2}+c'\|\nabla u^{n}\|^{2},\quad\forall n\ge 0,$$
 where $c'=c'(R_2)=L_g^2c_S^2+c(R_2)$. We multiply this by $2n\delta t$ and we add $\|\nabla u^{n+1}\|^2$ on both sides. This yields
$$(n+1)\|\nabla u^{n+1}\|^{2}\leq (1+2c'\delta t)n\|\nabla u^{n}\|^{2}+\|\nabla u^{n+1}\|^2,\quad \forall n\geq 0.$$
Let $d_{n}=n\|\nabla u^{n}\|^{2}$. By~\eqref{sobo},  we have
$$d_{n+1}\leq (1+2c'\delta t)d_{n}+ c_S\|\Delta u^{n+1}\|^{2},\quad \forall n\geq 0.$$
Using $d_{0}=0$, we deduce by induction that
$$d_{n}\leq (1+2c'\delta t)^{n}~(c_{S}\sum_{k=0}^{n-1}\|\Delta u^{k+1}\|^{2}),\quad \forall n\geq 1.$$
The conclusion~\eqref{c6}  follows from~\eqref{c2aux}  and from~\eqref{c21}  with $s=2c'\delta t$.
\end{proof}
\section{Finite time uniform error estimate}
\label{sec4}
For the error estimate on a finite time interval, we follow the methodology in \cite{P18,W10}.
We consider a sequence $(u^{n})$ in $D(A)$ generated by~\eqref{b1}. To the sequence $(u^{n})$, we associate three functions $u_{\delta t}$, $\overline{u}_{\delta t}$, $\underline{u}_{\delta t}$ : $\mathbb{R}_{+}\rightarrow D(A)$, namely 
\begin{eqnarray*}
&&u_{\delta t}= u^{n}+ \frac{t-n\delta t}{\delta t}(u^{n+1}-u^{n}),\quad  t\in [n\delta t,(n+1)\delta t),\\
&&\overline{u}_{\delta t}=u^{n+1},\quad t\in [n\delta t,(n+1)\delta t),\\
&&\underline{u}_{\delta t}= u^{n},\quad t\in [n\delta t,(n+1)\delta t).
\end{eqnarray*}
We note that $u_{\delta t}\in C^0([0,T],D(A))$ is piecewise linear, $ \overline{u}_{\delta t}\in L^\infty(0,T;D(A))$ and $\underline{u}_{\delta t}\in L^\infty(0,T;D(A))$, for all $T>0$. 
The scheme~\eqref{b1} can be rewritten
\begin{eqnarray}\label{d1}
\frac{du_{\delta t}}{dt}+A^{2}\overline{u}_{\delta t}+Af(\underline{u}_{\delta t})+ g(\underline{u}_{\delta t})=0\quad\mbox{ in }D(A^{-1}),\mbox{ for a.e. }t>0.
\end{eqnarray}
Equivalently, we have
\begin{eqnarray}
&&\frac{du_{\delta t}}{dt}+A^{2}u_{\delta t}+A f(u_{\delta t})+ g(u_{\delta t})\nonumber\\
&&=A^2(u_{\delta t}-\overline{u}_{\delta t}) + A(f(u_{\delta t})-f(\underline{u}_{\delta t}))+(g(u_{\delta t})-g(\underline{u}_{\delta t})) \label{d2}
\end{eqnarray}
in $D(A^{-1})$, for a.e. $t>0$. 
We denote by $u$ the solution to~\eqref{a1bis} with initial condition $u_0\in D(A)$ and we set $$e_{\delta t}(t)=u_{\delta t}(t)-u(t).$$
The error estimate reads:
\begin{theo}
	For all $T>0$ and for all $R_{2}>0$, there is a constant $C(T,R_{2})$ independent of $\delta t$ such that $u^{0}=u_{0}$ and $\|u^{0}\|_2\leq R_{2}$ imply
	\begin{equation}
	\sup_{t\in[0,N\delta t]}\|e_{\delta t}(t)\|\leq C(T,R_{2})~(\delta t)^{\frac{1}{2}},
	\end{equation}
		where $N=\lfloor T/\delta t \rfloor$ and $\lfloor \cdot \rfloor$ denotes the integer floor function.
\end{theo}
\begin{proof}
On subtracting~\eqref{a1bis} from~\eqref{d2}, we find
\begin{eqnarray}
&&\frac{de_{\delta t}}{dt}+\Delta^{2}e_{\delta t}+A(f(u_{\delta t})-f(u))+ (g(u_{\delta t})-g(u))\nonumber\\
&&=A^2(u_{\delta t}-\overline{u}_{\delta t}) + A(f(u_{\delta t})-f(\underline{u}_{\delta t}))+(g(u_{\delta t})-g(\underline{u}_{\delta t})) .\label{d3}
\end{eqnarray}
We multiply~\eqref{d3} by $e_{\delta t}$ in $H$. We obtain
\begin{eqnarray}
&&\frac {1}{2}\frac{d}{dt}\|e_{\delta t}(t)\|^{2}+ \|\Delta e_{\delta t}(t)\|^{2}-(f(u_{\delta t})-f(u),\Delta e_{\delta t})+(g(u_{\delta t})-g(u), e_{\delta t})\nonumber\\
&&=(\Delta(u_{\delta t}-\overline{u}_{\delta t}),\Delta e_{\delta t})-(f(u_{\delta t})-f(\underline{u}_{\delta t}),\Delta e_{\delta t})+(g(u_{\delta t})-g(\underline{u}_{\delta t}),e_{\delta t}).\label{d4}
\end{eqnarray}
Estimate~\eqref{c12bis} and Young's inequality yield
\begin{align*}
|(g(u_{\delta t})-g(u), e_{\delta t})|&\leq\|g(u_{\delta t})-g(u)\|~\| e_{\delta t}\|\\
&\leq L_{g}\| e_{\delta t}\|^{2}
\end{align*}
and
\begin{align*}
|(g(u_{\delta t})-g(\underline{u}_{\delta t}),e_{\delta t})|&\leq \|g(u_{\delta t})-g(\underline{u}_{\delta t})\|~\|e_{\delta t}\|\\
&\leq L_{g}\|u_{\delta t}-\underline{u}_{\delta t}\|~\|e_{\delta t}\|\\
&\leq \frac{L^{2}_{g}}{4}\|u_{\delta t}-\underline{u}_{\delta t}\|^{2}+\|e_{\delta t}\|^{2}.
\end{align*}
Moreover, by~\eqref{c12f},
\begin{align*}
|(f(u_{\delta t})-f(u),\Delta e_{\delta t})|&\leq \|f(u_{\delta t})-f(u)\|~\|\Delta e_{\delta t}\|\\
&\leq L_{f}\|e_{\delta}\|~\|\Delta e_{\delta t}\|\\
&\leq L^{2}_{f}\|e_{\delta}\|^{2}+\frac{1}{4}\|\Delta e_{\delta}\|^{2}
\end{align*}
and
\begin{align*}
|(f(u_{\delta t})-f(\underline{u}_{\delta t}),\Delta e_{\delta t})|&\leq \|f(u_{\delta t})-f(\underline{u}_{\delta t})\|~\|\Delta e_{\delta t}\|\\
&\leq L_{f}\|u_{\delta t}-\underline{u}_{\delta t}\|~\|\Delta e_{\delta t}\|\\
&\leq L^{2}_{f}\|u_{\delta t}-\underline{u}_{\delta t}\|^{2}+ \frac{1}{4}\|\Delta e_{\delta t}\|^{2}.
\end{align*}
We also have
\begin{align*}
|(\Delta(u_{\delta t}-\overline{u}_{\delta t}),\Delta e_{\delta t})|&\leq \|\Delta(u_{\delta t}-\overline{u}_{\delta t})\|~\|\Delta e_{\delta t}\|\\
&\leq \|\Delta(u_{\delta t}-\overline{u}_{\delta t})\|^{2}+\frac{1}{4}\|\Delta e_{\delta t}\|^{2}.
\end{align*}
Inserting the estimates above into~\eqref{d4}, we find
\begin{eqnarray}\label{d5}
\frac{d}{dt}\|e_{\delta t}(t)\|^{2}+ \frac{1}{2}\|\Delta e_{\delta t}(t)\|^{2}\leq c_{1}\|e_{\delta t}(t)\|^{2}+c_{2}\|u_{\delta t}-\underline{u}_{\delta t}\|^{2}+ 2\|\Delta(u_{\delta t}-\overline{u}_{\delta t})\|^{2},
\end{eqnarray}
where $c_{1}=2L_{g}+2+2L^{2}_{f}$ and $c_{2}=L_g^2/2+2L_f^2$.

Let $T>0$ and $N=\lfloor T/\delta t \rfloor$. Thanks to $e_{\delta t}(0)=0$ and the classical Gronwall lemma applied to~\eqref{d5}, we obtain
\begin{eqnarray}
\|e_{\delta t}(t)\|^{2}&\leq& \exp{(c_{1}T)}\int_{0}^{N\delta t}c_{2}\|u_{\delta t}(s)-\underline{u}_{\delta t}(s)\|^{2}ds\nonumber\\
&& + \exp{(c_{1}T)}\int_{0}^{N\delta t}2\|\Delta(u_{\delta t}(s)-\overline{u}_{\delta t}(s))\|^{2}ds,\quad \forall t\in[0, N\delta t].\label{d6}
\end{eqnarray}
 On the interval $[n\delta t, (n+1)\delta t)$, we have
$$\|u_{\delta t}(s)-\underline{u}_{\delta t}(s)\|\leq \|u^{n+1}-u^{n}\|\quad\mbox{ and }\quad \|\Delta(u_{\delta t}(s)-\overline{u}_{\delta t}(s))\|\leq \|\Delta(u^{n+1}-u^{n})\|.$$
 Thus,
$$\int_{0}^{N\delta t}c_{2}\|u_{\delta t}(s)-\underline{u}_{\delta t}(s)\|^{2}ds\leq c_{2}\delta t\sum_{k=0}^{N-1}\|u^{k+1}-u^{k}\|^{2}$$
and
$$\int_{0}^{N\delta t}\|\Delta(u_{\delta t}(s)-\overline{u}_{\delta t}(s))\|^{2}ds\leq \delta t\sum_{k=0}^{N-1}\|\Delta(u^{k+1}-u^{k})\|^{2}.$$
Plugging these estimates into~\eqref{d6}, we obtain
\begin{eqnarray}\label{d7}
\|e_{\delta t}(t)\|^{2}&\leq& \exp{(c_{1}T)}\left(c_{2}\sum_{k=0}^{N-1}\|u^{k+1}-u^{k}\|^{2}+2\sum_{k=0}^{N-1}\|\Delta(u^{k+1}-u^{k})\|^{2}\right)\delta t\nonumber\\
&\le &c_{3}\exp{(c_{1}T)}\sum_{k=0}^{N-1}\|\Delta(u^{k+1}-u^{k})\|^{2}\delta t,
\end{eqnarray}
where $c_{3}=c_{2}c_S^2+2$. By~\eqref{dissiH2dtbis}, 
$$\|e_{\delta t}(t)\|^{2}\le c_{3}\exp{(c_{1}T)}\left(Q_2(\|\Delta u^0\|)+M_2'T\right)\delta t,\quad \forall t\in[0, N\delta t].$$
This concludes the proof. 
\end{proof}
\section{Convergence of exponential attractors}
\subsection{Some  definitions}
Before stating our main result, we recall some  definitions (see e.g.~\cite{EMZ,T}).  We recall that $H=L^2(\Omega)$ and $D(A)=H^2(\Omega)\cap H^1_0(\Omega)$.
  A continuous-in-time semigroup $\{S(t),\ t\in\Rr_+\}$ on  $D(A)$ is a family of (nonlinear) operators such that $S(t)$ is a continuous operator (for the $L^2(\Omega)$-norm) from $D(A)$ into itself, for all $t\in\Rr_+$, with $S(0)=Id$ (identity)     and  
$$S(t+s)=S(t)\circ S(s),\quad\forall s,t\in\Rr_+.$$
A discrete-in-time semigroup $\{S(t),\ t\in\Nn\}$ on $D(A)$ is a family of (nonlinear) operators which satisfy these properties with $\Rr_+$ replaced by $\Nn$.
A discrete-in-time semigroup is usually  denoted $\{S^n,\ n\in\Nn\}$, where $S(=S(1))$ is a  continuous (nonlinear) operator from $D(A)$ into itself.  

A (continuous or discrete) semigroup $\{S(t),\ t\ge 0\}$ defines a (continuous or discrete) dynamical system: if $u_0$ is the state of the dynamical system at time $0$, then $u(t)=S(t)u_0$ is the state at time $t\ge 0$.  The term ``dynamical system'' will sometimes be used instead of ``semigroup''.

\begin{definition}[Global attractor] Let $\{S(t),\ t\ge 0\}$ be a continuous or discrete semigroup on  $D(A)$. A bounded set $\mathcal{A}\subset D(A)$ is called the global attractor of the dynamical system if the following three conditions are satisfied:
\begin{enumerate}
\item $\mathcal{A}$ is compact in $H$;
\item $\mathcal{A}$ is invariant, i.e. $S(t)\mathcal{A}=\mathcal{A}$, for all $t\ge 0$;
\item $\mathcal{A}$ attracts all bounded sets in $D(A)$, i.e., for every bounded set $B$ in $D(A)$, 
$$\lim_{t\to+\infty} dist_H(S(t)B,\mathcal{A})=0.$$
\end{enumerate}
Here, $dist_H$ denotes the non-symmetric Hausdorff semidistance in $H$ between two subsets, which is defined as 
$$dist_H(B_1,B_2)=\sup_{b_1\in B_1}\inf_{b_2\in B_2}\|b_1-b_2\|_H.$$
\end{definition}
It is easy to see, thanks to the invariance and the attracting property, that the global attractor, when it exists, is unique~\cite{T}. 

Let $X\subset H$ be a (relatively compact) subset of $H$. For $\ep>0$, we denote $N_\ep(X,H)$ the minimum number of balls of $H$ of radius $\ep>0$ which are necessary to cover $X$. The {\it fractal dimension of }$X$ (see e.g.~\cite{EFNT,T}) is the number
$$dim_F(X)=\limsup_{\ep\to 0}\frac{\log(N_\ep(X,H))}{\log(1/\ep)}\in[0,+\infty].$$

\begin{definition}[Exponential attractor] Let $\{S(t),\ t\ge 0\}$ be a continuous or discrete semigroup on  $D(A)$. A bounded set $\mathcal{M}\subset D(A)$ is an exponential attractor of the dynamical system if the following three conditions are satisfied:
\begin{enumerate}
\item $\mathcal{M}$ is compact in $H$ and has finite fractal dimension;
\item $\mathcal{M}$ is positively invariant, i.e. $S(t)\mathcal{M}\subset\mathcal{M}$, for all $t\ge 0$;
\item $\mathcal{M}$ attracts exponentially the bounded subsets of $D(A)$ in the following sense:
$$\forall B\subset D(A)\mbox{ bounded},\quad dist_H(S(t)B,\mathcal{M})\le\mathcal{Q}(\|B\|_H)  e^{-\alpha t},\quad t\ge 0,$$
where the positive constant $\alpha$ and the monotonic function $\mathcal{Q}$ are  independent of $B$. Here, $\|B\|_H=\sup_{b\in B}\|b\|_H$.
\end{enumerate}
\end{definition}
It is easy to see that the exponential attractor, if it exists, contains the global attractor.

\subsection{The main result}
We have seen  that $\{S_0(t),\ t\in\Rr_+\}$ defined by~\eqref{defSt} is  a continuous-in-time dynamical system on $D(A)$, and that for every ${\delta t}>0$ small enough, $\{S^n_{\delta t},\ n\in\Nn\}$ defines a discrete-in-time dynamical system on $D(A)$ (Theorem~\ref{b2}). We have:
\begin{theo}
\label{theomain}
Let $\delta t^\star>0$ be small enough.  For every ${\delta t}\in (0,{\delta t}^\star]$, the discrete dynamical system $\{S_{\delta t}^n,\ n\in\Nn\}$ possesses an exponential attractor $\mathcal{M}_{\delta t}$ in $D(A)$, and the continuous dynamical system $\{S_0(t),\ t\in\Rr_+\}$ possesses an exponential attractor $\mathcal{M}_0$ in $D(A)$ such that:
\begin{enumerate}
\item the fractal dimension of $\mathcal{M}_{\delta t}$ is bounded, uniformly with respect to ${\delta t}\in[0,{\delta t}^\star]$, 
$$dim_F\mathcal{M}_{\delta t}\le c_{1},$$
where $c_{1}$ is independent of ${\delta t}$;
\item $\mathcal{M}_{\delta t}$ attracts the bounded sets of $D(A)$, uniformly with respect to ${\delta t}\in(0,{\delta t}^\star]$, i.e. for all ${\delta t}\in(0,{\delta t}^\star]$,
$$\forall B\subset D(A)\mbox{ bounded },\ dist_H(S^n_{\delta t} B,\mathcal{M}_{\delta t})\le \mathcal{Q}(\|B\|_H)e^{-c_{2}n{\delta t}},\quad n\in\Nn,$$
where the positive constant $c_{2}$ and the monotonic function $\mathcal{Q}$ are independent of ${\delta t}$; 
\item the family $\{\mathcal{M}_{\delta t},\ {\delta t}\in[0,{\delta t}^\star]\}$ is continuous at $0$,
$$dist_{sym}(\mathcal{M}_{\delta t},\mathcal{M}_0)\le c_{3}(\delta t)^{c_{4}},$$
where $c_{3}$ and $c_{4}\in(0,1)$ are independent of ${\delta t}$ and $dist_{sym}$ denotes the symmetric Hausdorff distance between sets, defined by
$$dist_{sym}(B_1,B_2):=\max\{dist_H(B_1,B_2),dist_H(B_2,B_1)\}.$$
\end{enumerate}
\end{theo}
\begin{proof}We apply Theorem 2.5 in~\cite{P18} with the spaces $H=L^2(\Omega)$ and $V=H^1_0(\Omega)$ and the set 
$$\mathcal{B}=\left\{ v\in D(A)\ :\ \|v\|_2\le M_2+1\right\},$$
where $M_2$ is the constant in~\eqref{a27bis} and~\eqref{dissiH2dt}. 
  We note that $V$ is compactly imbedded in $H$ and that an $H$-$V$ smoothing property holds, uniformly with respect to $\delta t$ (Lemma~\ref{a35} and Lemma~\ref{lem30}). Moreover, $\mathcal{B}$  is absorbing in $D(A)$, uniformly with respect to ${\delta t}\in[0,{\delta t}_0]$, where $\delta t_0>0$ is chosen small enough.  The estimates of Sections~\ref{sec2}-\ref{sec4} show that assumptions (H1)-(H9) of Theorem 2.5 in~\cite{P18} are satisfied.  Thus, the conclusions of Theorem~\ref{theomain} hold for ${\delta t}\in[0,{\delta t}^\star]$, for some ${\delta t}^\star\in(0,{\delta t}_0]$  small enough.  We note that Theorem 2.5 in~\cite{P18} is stated for a family of semigroups which act on the whole space $H$, but with a minor modification of the proof, it can be applied  to our situation where the semigroup acts on $D(A)$ and is continuous for the $H$-norm. The main tool is the construction of exponential attractors based on a uniform smoothing property proposed by  Efendiev, Miranville and Zelik  in~\cite[Theorem 4.4]{EMZ}. 
\end{proof}
As in~\cite[Corollary 6.2]{P18}, we have:
\begin{coro}For every ${\delta t}\in[0,{\delta t}^\star]$, the semigroup $\{S_{\delta t}(t),\ t\ge 0\}$ possesses a global attractor $\mathcal{A}_{\delta t}$ in $D(A)$ which is bounded in $D(A)$ and  compact  in $H$. Moreover, $dist_H(\mathcal{A}_{\delta t},\mathcal{A}_0)\to 0$ as ${\delta t}\to 0^+$, and the fractal dimension of $\mathcal{A}_{\delta t}$ is bounded by a constant independent of ${\delta t}$. 
\end{coro}

\begin{rem}
\label{remNeu}
Let us  replace the Dirichlet boundary conditions~\eqref{a2} with Neumann boundary conditions, which read
\begin{equation}
\label{a2bis}
\partial_n u=\partial_n\Delta u=0\mbox{ on }\partial \Omega\times (0,+\infty),
\end{equation}
where $n$ is the unit outer normal to $\partial\Omega$. 

 If $g=0$ in~\eqref{a1}, we deal with the classical Cahn-Hilliard equation and a result similar to Theorem~\ref{theomain} can be obtained. In this case, we have the conservation of  mass and it is convenient to introduce the function spaces 
 $$H_\beta=\left\{ v\in H\ :\ \int_\Omega v=\beta\right\}\mbox{ and }\mathcal{H}_\alpha=\bigcup_{|\beta|\le \alpha}H_\beta,$$
 as in~\cite{BP18,T}, where $\beta\in\Rr$ and $\alpha>0$.

If $g\not=0$, the situation is more delicate because we no longer have the conservation of mass~\cite[Remark 5.7]{M13}. If $g$ is a proliferation term, the mass may even blow in finite time~\cite{CMZ14,MAIMS,Mbook}. 
\end{rem}

\end{document}